\documentclass[smallextended]{svjour3}
\smartqed
\usepackage[margin=1in]{geometry}
\usepackage[onehalfspacing]{setspace}
\usepackage[utf8]{inputenc}

\usepackage[english]{babel}

\usepackage{pdflscape}
\usepackage{fancyvrb}
\usepackage{calc}
\usepackage{rotating}
\usepackage{pgf,tikz}
\usepackage{mathrsfs}
\usetikzlibrary{arrows}
\usepackage{mathtools}
\DeclarePairedDelimiter\ceil{\lceil}{\rceil}

\usepackage{listings}
\usepackage[]{qcircuit}
\usepackage{braket}
\usepackage{mathtools}
\usepackage{varwidth}
\usepackage{caption}
\usepackage{subcaption}
\usepackage{graphicx}

\makeatletter
\def\paragraph{\@startsection{paragraph}{4}%
  \z@\z@{-\fontdimen2\font}%
  {\normalfont\bfseries}}
\makeatother

\usepackage{graphicx}
\usepackage{pgffor}
\usepackage{tikz-cd}
\usepackage{booktabs}
\usepackage{enumerate}

\usepackage{amssymb, amsmath, amsthm}
\usepackage[linesnumbered,ruled]{algorithm2e}

\usepackage{scalerel}
\usepackage{ulem}

\usepackage{graphicx}              
\usepackage{amsmath}               
\usepackage{amsfonts}              
\usepackage{amsthm}                
\usepackage{xkeyval}               
\usepackage{amssymb}
\usepackage{enumerate}             
\usepackage{color}
\usepackage{breqn}                 
\usepackage{bm}                    

\DeclareUnicodeCharacter{24C7}{}

\allowdisplaybreaks[1]
\usepackage{nicefrac}
\usepackage{booktabs}

\usepackage{kpfonts}
\usepackage{stackengine}
\usepackage{calc}
\newlength\shlength
\newcommand\xshlongvec[2][0]{\setlength\shlength{#1pt}%
	\stackengine{-5.6pt}{$#2$}{\smash{$\kern\shlength%
			\stackengine{7.55pt}{$\mathchar"017E$}%
			{\rule{\widthof{$#2$}}{.57pt}\kern.4pt}{O}{r}{F}{F}{L}\kern-\shlength$}}%
	{O}{c}{F}{T}{S}}

\usepackage{hyperref}
\hypersetup{
	colorlinks   = true,    
	citecolor    = red      
}

\newcommand{\RN}[1]{%
  \textup{\uppercase\expandafter{\romannumeral#1}}%
}

\allowdisplaybreaks

\newcommand{\meqref}[1]{Eq.~$\! \eqref{#1} $}
\newcommand{\mref}[1]{Sect.~$ \!\ref{#1} $}
\newcommand{\mfig}[1]{Fig.~$ \!\ref{#1} $}

\newcommand{\RR}{\mathbb{R}}      

\newcommand{\mat}[4]{\left[\begin{smallmatrix*}[r]
                #1 & #2 \\
                #3 & #4 \\
        \end{smallmatrix*}\right]}

\def\CC{{\mathbb C}}

\def\RR{{\mathbb R}}
\def\SS{{\mathbb S}}
\def\TT{{\mathbb T}}

\def\ll{{\frak l}}  

\def\ss{{\frak s}}
\def\tt{{\frak t}}

\newcommand{\myN}{\zeta}

\def\<{\langle}
\def\>{\rangle}

\def\sl{\ss\ll}

\numberwithin{equation}{section}

\begin{document}

	\title{Trajectory optimization using quantum computing}
	
	\author{Alok Shukla       \and
		Prakash Vedula 
	}
	
	
	\institute{Alok Shukla  \at
		The University of Oklahoma, Norman, USA \\
		\email{alok.shukla@ou.edu}           
		\and
		Prakash Vedula \at
	   The University of Oklahoma, Norman, USA \\
		\email{pvedula@ou.edu}           
}

	\date{}

	\maketitle
	
	\begin{abstract}

We present a framework wherein the trajectory optimization problem (or a problem involving calculus of variations) is formulated as a search problem in a discrete space.  A distinctive feature of our work is the treatment of discretization of the optimization problem wherein we discretize not only independent variables (such as time) but also dependent variables. Our discretization scheme enables a reduction in computational cost through selection of coarse-grained states. It further facilitates the solution of the trajectory optimization problem via classical discrete search algorithms including deterministic and stochastic methods for obtaining a global optimum. This framework also allows us to efficiently use quantum computational algorithms for global trajectory optimization.  We demonstrate that the discrete search problem can be solved by a variety of techniques including a deterministic exhaustive search in the physical space or the coefficient space, a randomized search algorithm, a quantum search algorithm or by employing a combination of randomized and quantum search algorithms depending on the nature of the problem. We illustrate our methods by solving some canonical problems in trajectory optimization. We also present a comparative study of the performances of different methods in solving our example problems. Finally, we make a case for using quantum search algorithms as they offer a quadratic speed-up in comparison to the traditional non-quantum algorithms. 

	\end{abstract}

\section{Introduction}

The goal of trajectory optimization 
is to find a path or trajectory that optimizes a given quantity of interest or any other objective function associated with a certain performance measure, under a set of given constraints on the dynamics of the system. Trajectory optimization problems appear naturally in many practical situations, especially in aerospace applications. Trajectory optimization problems are important, even at a more fundamental level, as they can be related to the principle of least action or other variational principles underlying diverse physical phenomena (see \cite{feynman2005principle}, \cite{goldstein2011classical}, \cite{gelfand2000calculus}, \cite{elsgolc2012calculus}).

Mathematically, action is a functional with the trajectory function of the system as its argument. Minimizing the action functional means selecting the optimum trajectory that minimizes the underlying integral associated to the action. Of course, in this sense, trajectory optimization is related to the subject of the calculus of variations. Therefore, the techniques proposed in this paper are also applicable to problems involving the calculus of variations.

There are many known methods of trajectory optimization (see \cite{betts1998survey}, \cite{cowling2006optimal}, \cite{fahroo2002direct}, \cite{von1992direct}). These solution methods can be broadly classified into the categories of \texttt{indirect methods} and  \texttt{direct methods}. In an {indirect method} one proceeds by analytically finding a set of necessary and sufficient conditions to solve a given trajectory optimization problem. {Indirect methods} have a long history dating back to the celebrated brachistochrone problem, posed by Bernoulli in the $ 16^{\text{th}}$ century (\cite{bernoulli1996problema}, \cite{woodhouse1810treatise}), and which was subsequently solved by Euler, Lagrange and others by employing the techniques of calculus of variations. On the other hand {direct methods} have become more popular with the advent of digital computers. A  {direct method} consists of discretizing the optimization problem and then solving the resulting non-linear optimization problem directly. Typically, in a {direct method} the control and the state parameters are represented by piecewise polynomial functions or a linear combination of global basis functions satisfying the boundary constraint on a discrete time grid. Several types of polynomials have been employed to represent the control and state parameters in the literature. For example in \cite{fahroo2002direct}, Chebyshev polynomials are used to solve the optimization problem. 

We remark that  {indirect methods} often involve  prior mathematical analysis of the problem. The resulting analytical solution is typically quite complicated and not often amenable to a computationally efficient solution.  Therefore, in practical applications {direct methods} are often more suitable. However, in high dimensional problems direct methods also become computationally expensive. In this paper we describe a novel discretization scheme to solve the global optimization problem wherein we discretize  both the independent variables (such as time)  and the dependent variables.  Our discretization scheme involves the selection of a finite number of coarse-grained states.  Each of these states has an associated cost function. The problem then becomes a discrete search problem, in which a state with the minimum associated cost is to be determined.  Then it becomes possible to employ deterministic and probabilistic methods for obtaining a global optimum. In this work, we propose a new framework for solution of the trajectory optimization problem via classical discrete search algorithms including an exhaustive search algorithm (Method I, \mref{subsect_classical_exhaustive}), a random search algorithm (Method II, \mref{random}), and a hybrid search algorithm (Method III, \mref{subsect_classical_hybrid_search}). This framework also allows us to efficiently use quantum computational algorithms for global trajectory optimization.  In this context, we propose new approaches for solution of the trajectory optimization problem using quantum exhaustive search algorithms (Method IV, \mref{subsect5_quantum_exhaustive}), a quantum random search algorithm (Method V, \mref{subsect5_quantum_random}), and a quantum hybrid algorithm (Method VI, \mref{subsect5_quantum_hybrid}). It turns out that quantum computers, in principle, are significantly superior to classical computers in solving the underlying discrete search problems. In fact, a slight modification of Grover's quantum search algorithm (see \cite{grover1996fast}) leads to a solution with complexity of the order of $ O(\sqrt{N}) $, where $ N $ depends on the number of possible states of control and state parameters. The number $ N $ is more precisely defined in \mref{Discretization}.

A main focus of this paper is to show that trajectory optimization problems can be tackled efficiently using quantum computational algorithms employed either alone or in conjunction with a randomized search. As noted earlier, to achieve this we use a discretization scheme which makes the use of quantum computing possible. We also note here that unlike many other works in the literature on trajectory optimization which are set up to only detect the local optimum, our approach enables the search of the global optimum. We will demonstrate our method using two canonical problems in trajectory optimization, namely the brachistochrone problem and the moon landing problem.  The method presented here could be made even more effective if it is combined with other techniques like the gradient descent and simulated annealing (see \cite{vanderbilt1984monte}). However, those aspects are outside the scope of the present work and will be considered in a future work.

In \mref{Problem formulation}, we present the problem formulation and in the following section, \mref{Discretization}, we describe our discretization scheme. In \mref{sect_classical_search_algo} we describe classical search algorithms including the exhaustive, random and hybrid search algorithms to solve the discrete version of the trajectory optimization problem formulated in \mref{Discretization}. \mref{Quantum}  contains a brief introduction to quantum search methods, including Grover's algorithm. We also describe various quantum search algorithms in this section and discuss advantages of these algorithms over their classical counterparts. \mref{Computational examples}  contains three computational examples, namely, the brachistochrone problem, the isoperimetric problem and the moon landing problem. In this section, we give several approaches to discretization and searching such as discretization in the physical space or the co-efficient space, followed by an exhaustive search, a random search, a quantum search or a hybrid search. We also discuss and compare the performances of these approaches in solving the brachistochrone problem. Finally, we present our concluding remarks in \mref{Conclusion}.

\section{Trajectory optimization: Problem formulation} \label{Problem formulation}
We will closely follow \cite{fahroo2002direct} in the formulation of the trajectory optimization problem in this paper. Suppose $ \vec{U}(t) \in \RR^m  $ and $ \vec{X}(t) \in \RR^n  $ denote the control function  and the corresponding state trajectory respectively at time $ t $. The goal is to determine the optimal control function $ \vec{U}(t)$ and the corresponding state trajectory $ \vec{X}(t)$ for $ \tau_0 \leq t \leq \tau_f $ such that the following Bolza cost function is minimized:
\begin{align}\label{Eq_Main_Cost}
\mathcal{J}(\vec{U}(.),\vec{X}(.),\tau_f) = \mathcal{M}(\vec{X}(\tau_f),\tau_f) + \int_{\tau_0}^{\tau_f} \, \mathcal{L}(\vec{U}(\tau),\vec{X}(t),t) \, dt.
\end{align}
Here $ \mathcal{M} $ and $ \mathcal{L} $ are $ \RR $-valued functions. Moreover, we also assume that the system satisfies the following dynamic constraints. 
\begin{align}
f_l \leq f(\vec{U}(\tau),\vec{X}(t),\vec{X^'}(t),t) \leq f_u \qquad t \in [\tau_0, \tau_f].
\end{align}
In addition, we also specify the following boundary conditions 
\begin{align}
h_l  \leq h(\vec{X}(\tau_0),\vec{X}(\tau_f),\tau_f -\tau_0) \leq h_u,
\end{align}
where $ h $ an $ \RR^p $-valued function and $ h_l, h_u \in \RR^p$ are constant vectors providing the lower and upper bounds of $ h $. Finally, we note the mixed constraints on control and state variables
\begin{align}
g_l \leq g(\vec{U}(t),\vec{X}(t),t) \leq g_u ,
\end{align}
with $ g $ a $ \RR^r $-valued function and $ g_l, g_u \in \RR^r$ are constant vectors providing the lower and upper bounds of $ g $.

We note that with an appropriate transformation we may assume that $ \tau_0 =-1 $ and $ \tau_f = 1$ (see \cite{fahroo2002direct}). 

\section{Discretization approaches} \label{Discretization}
There are a number of local and global discretization methods available for the problem described in the previous section. Often the trajectory is approximated using a linear combination of a set of orthogonal polynomials. For example, one can use global spectral methods employing Chebyshev polynomials. Another possibility is to use piecewise quadratic or cubic splines to approximate the trajectory within a set of chosen discrete grid points. An appropriate choice of discretization method also depends upon the nature of the problem. For example, the Chebyshev pseudospectral method has the advantage that the node points cluster around the end points of the interval avoiding the Runge phenomenon. In this paper we present a piecewise formulation which has some advantages over global polynomial approximations. Usually, one finds in the literature that the discretization is carried out in the independent variable such as time (or space). In addition to discretizing the independent variable, in this paper we propose to discretize the dependent variable as well. In other words, discretization will mean mapping a set of continuous dependent (and independent) variables to a set of discrete dependent (and independent) variables. Next we describe our discretization scheme in detail.

We divide the time interval $ [\tau_0,\tau_f] $ into $ \eta $ sub-intervals $ [t_j,t_{j+1}] $ for $ j=0,1 \cdots \eta-1  $ with $ \tau_0 = t_0 < t_1 < t_2 < \cdots < t_{\eta-1}  < t_\eta =\tau_f $.
We define 
\begin{align}
\vec{X}(t) = \sum_{j=0}^{\eta-1} \,  \boldsymbol{\chi}_j(t) 
\end{align}
with 
\begin{align}
\boldsymbol{\chi}_j(t) = \begin{cases} \displaystyle 
\sum_{k=0}^{m_j} \, \vec{a}_{j,k} \, \phi_{j,k}(t)  \qquad &\text{if }  t \in [t_j,t_{j+1}), \\ \\
0   &\text{otherwise}.
\end{cases}
\end{align}
Here $ \vec{a}_{j,k} \in \RR^n $ are constant vectors and $ \phi_{j,k} $ is a chosen spectral polynomial for approximating the trajectory within $ t \in [t_j,t_{j+1}] $. For example, for a quadratic approximation we may assume $ m_j =2 $ and $ \phi_{j,k} = t^k $. Further, depending upon the nature of the problem, the following conditions may be needed to avoid any possible discontinuities
\begin{align} \label{eq_boundary_discrete}
\lim\limits_{t \to  t_{j+1}} \boldsymbol{\chi}_j(t) = \boldsymbol{\chi}_{j+1}(t_{j+1}) \qquad \text{for } j=0,1,2 \ldots \eta-1.
\end{align}
Moreover, if needed, smoothness conditions involving higher-order derivatives may also be imposed. 

We make a similar definition for $ \vec{U} $ as,
\begin{align}
\vec{U}(t) = \sum_{j=0}^{\eta-1} \,  \boldsymbol{\mu}_j(t) 
\end{align}
with 
\begin{align}
\boldsymbol{\mu}_j(t) = \begin{cases} \displaystyle 
\sum_{k=0}^{n_j} \, \vec{b}_{j,k} \, \psi_{j,k}(t)  \qquad &\text{if }  t \in [t_j,t_{j+1}), \\ \\
0   &\text{otherwise}.
\end{cases}
\end{align}
Here $ \vec{b}_{j,k} \in \RR^m $ are constant vectors and $ \psi_{j,k} $ is a chosen spectral polynomial for approximating the control function within $ t \in [t_j,t_{j+1}] $. Further, we impose the following boundary conditions similar to  \meqref{eq_boundary_discrete} 
\begin{align} \label{eq_boundary_discrete1}
\lim\limits_{ t \to t_{j+1}}\boldsymbol{\mu}_j(t) = \boldsymbol{\mu}_{j+1}(t_{j+1}) \qquad \text{for } j=0,1,2 \ldots \eta-1.
\end{align}
We note that the above formulation is more general than the one given in \cite{fahroo2002direct} as on setting  $ m_j =0 $ and suitably defining $ \chi $ and $ \mu $ we can recover the formulation given in \cite{fahroo2002direct}.

Now we further discretize our problem by imposing the condition that for $ 1 \leq r \leq n $ and for $ 0 \leq k \leq m_j $  the components of vectors $ \vec{a}_{j,k} (r) $  take values from a discrete set of cardinality $ S_{j,r,k} $, say $ \vec{a}_{j,k}(r) \in \SS_{j,r,k} $ where
$ \SS_{j,r,k} = \{ \alpha_1, \alpha_2, \cdots   \} $ with $  \alpha_1, \alpha_2, \cdots \in \RR $ and 
$ \#\SS_{j,r,k} = S_{j,r,k} $. We note that the choices of $ \alpha_i $, the elements of the set  $ \SS_{j,r,k} $ 
depend on the nature of the problem and can be appropriately modified. A typical choice could be equidistant entries of the following form
$$ \SS_{j,r,k} = \{ -n_{j,r,k},-(n_{j,r,k}-\epsilon_{j,r,k}), \ldots,0, \ldots,n_{j,r,k}-\epsilon_{j,r,k},n_{j,r,k}\} \subset \RR.$$
Similarly, for $ 1 \leq r \leq m $ and for $ 0 \leq k \leq n_j $ we let the components of vectors $ \vec{b}_{j,k} (r) $  to take values only from a discrete set of cardinality $ T_{j,r,k} $, say $ \vec{b}_{j,k}(r) \in \TT_{j,r,k} $, where 
$ \TT_{j,r,k} = \{ \beta_i \ | \ \beta_i \in \RR \} $
with $ \#\TT_{j,r,k} = T_{j,r,k} $.

We note that with this discretization the trajectory optimization problem is now turned into a discrete search problem of size $N = ST  $ where 
\begin{align} \label{Eq_S}
S =  \prod_{j=0}^{\eta-1}\prod_{r=1}^{n} \prod_{k=0}^{m_j} S_{j,k,r}  
\end{align} 
and
\begin{align} \label{Eq_T}
T =  \prod_{j=0}^{\eta-1}\prod_{r=1}^{m} \prod_{k=0}^{n_j} T_{j,k,r} \,, 
\end{align} 
and where the objective function to be minimized is the discrete form of  \meqref{Eq_Main_Cost}.
Once we take the boundary conditions given by  \meqref{eq_boundary_discrete}   and  \meqref{eq_boundary_discrete1}  into account, the problem size is further reduced. For example, for $ j=1 \ldots \eta-1$ we can always choose $ \vec{a}_{j,0} $ and $ \vec{b}_{j,0} $  such that the boundary conditions given by  \meqref{eq_boundary_discrete}   and  \meqref{eq_boundary_discrete1}  are satisfied. In any case it is clear that there are only a finite number of states, say $ N $, in our search space. 

In principle, we can now solve the problem by traversing our discrete search space and finding the state that corresponds to the minimum cost and satisfies all the imposed constraints. Still, in practice, the problem is the massive size of the discrete search space. We remark that even in the double precision $ 64 $-bit representation of a variable (a real number) on a digital computer, a total of $ 2^{64} $ distinct states of the variable can be faithfully represented on the machine. In an optimization problem involving hundreds of state variables, clearly the size of the discrete search space could become very large and unmanageable. We propose to use a discretization scheme consisting of ``coarse-grained''  states, thereby reducing the size of the discrete search space compared to the traditional ``fine-grained'' states  representation on a digital computer. The trade-off here is that reducing the size of the discrete search space affects the quality of the solution. Therefore, it becomes desirable to  obtain a delicate balance between the size of the discrete search space and the quality/accuracy of the resulting solution.  In most applications, the size of the discrete search space would render a direct exhaustive search computationally very costly or even outright impossible. In such cases we look for alternate practical solutions. One may use analytical techniques like the gradient method to appropriately direct the search.  Another possible approach is to use a random search method to find an approximate solution and then use an exhaustive search near that approximate solution. Parallel computing may also be helpful in this context. In this article, we propose to employ a quantum search algorithm which is faster and more efficient in comparison to all the classical methods. We also discuss the possibility of using a hybrid algorithm combining random search with the quantum search. But first we consider classical approaches to solve global trajectory optimization problems.   
\section{\textbf{New framework for trajectory optimization via  classical search algorithms}}\label{sect_classical_search_algo}
  In this section, we consider the problem described in \mref{Discretization} of finding the minimum cost, with the cost function being the discrete version of  \meqref{Eq_Main_Cost}  and $ N = ST $ with $ S $ and $ T $ as in  \meqref{Eq_S}  and  \meqref{Eq_T}.
  We will describe the application of classical algorithms, including deterministic and probabilistic approaches for solving this problem using our proposed framework of coarse-grained states. We remark that the optimum given by these methods are global, unlike many instances in the trajectory optimization literature wherein local optima are obtained.    
\subsection{\textbf{Method I: Classical exhaustive search algorithm}} \label{subsect_classical_exhaustive}
  We say that the search is exhaustive if each of the $ N $ states have been searched for the optimum cost. An exhaustive search in the physical space is guaranteed to give an optimal global solution to any desired degree of accuracy with appropriate grid refinement. In  \mref{sub_sect_discrete_physical_space}  we will present further details on this method in the context of solving the brachistochrone problem.

\subsection{\textbf{Method II: Classical random search algorithm}} \label{random}
 In many trajectory optimization problems the discrete search space is so large that an exhaustive search becomes computationally very costly. An alternative approach in such cases is to apply the well-known ``pure random search'' method, wherein samples are randomly selected and the running minimum cost is updated if needed (see \cite{brooks1958discussion}, \cite{patel1989pure}, \cite{zabinsky1992pure}).  In our context this algorithm is described below in Algorithm~$ \ref{Algo_one} $. In  \mref{subsect_random}, we will present a computational example employing this algorithm for solving the brachistochrone problem.
 
 \begin{algorithm}[ht]
 	\SetKwInOut{Input}{Input}
 	\SetKwInOut{Output}{Output}
 	\tcc{	Description: The objective of this random algorithm is find the minimum cost for the problem described in Sect.$ \! \ref{Discretization}$ with the cost function being the discrete version of  \meqref{Eq_Main_Cost}  and $ N = ST $ with $ S $ and $ T $ as in  \meqref{Eq_S}  and  \meqref{Eq_T}.}
 	\underline{function findMinCostRandom} ($ \{Y_0,Y_1,Y_2 \ldots ,Y_ {N-1}   \}, n $)\;
 	\Input{a list of size $ N $ containing all the possible states of the system and $ n $ the number of times the search has to be performed.}
 	\Output{the index of an admissible state which has the minimum cost in $ n $ random trials, the associated minimum cost.}
 	Set $ x = $  a random state from the input list.   \\
 	Set $ count = 1 $. \\
 	\While{$ count < n $}{
 		Set $ y = $  a random state from the input list. \\ 
 		Set $ count = count +1 $.\\
 		\If{$ cost(y) < cost(x) $.}
 		{
 			Set $ x = y $. \\
 		}
 	}
 	Return ($ x $, $ cost(x) $).
 	\caption{A random search algorithm for trajectory optimization} \label{Algo_one}
 \end{algorithm}

We note that this random search algorithm is guaranteed to succeed, with probability $ 1 $, as the number of steps $ n $ approaches infinity. To make this more precise, suppose that all the states are equally likely to yield the minimum cost. Let $ p(r,n) $ denote the probability of the random search algorithm picking the $ r^{\text{th}} $ ranked minimum cost path in the list, in $ n $ random searches. Clearly $  p(r,n) = \sum_{k=0}^{n-1} (\frac{r}{N})(\frac{N-r}{N})^k  = 1 - (1-\frac{r}{N})^n$ and  $ \lim\limits_{n \to \infty} p(r,n) \to 1 $.
 
 \subsection{\textbf{Method III: Classical hybrid search algorithm}} \label{subsect_classical_hybrid_search}  The hybrid search algorithm is a combination of the random and the exhaustive search algorithms and in many situations it is the most efficient classical search method. In a hybrid search algorithm first a random search is carried out to obtain an approximate solution. Then the result of the random search is further improved by performing an exhaustive search in a finer discrete grid near the approximate solution given by the random search. A computational example will be given in \mref{subsect_hybrid}.
 
\section{\textbf{New framework for trajectory optimization via quantum search algorithms}} \label{Quantum}
Quantum computing is an exciting field which beautifully combines quantum physics with computer science. There has been impressive progress in theoretical development of quantum algorithms as well as in their practical implementation on quantum computers. We refer the readers to any standard book on the subject for a more in-depth treatment (for eg., \cite{nielsen2002quantum}, \cite{rieffel2011quantum} or \cite{yanofsky2008quantum} ).  

\subsection{\textbf{Overview of quantum computation}}
A fundamental concept in quantum computation is that of the quantum bit, or `qubit'.  Unlike the classical bit, the qubit, say $ \ket{q}$, can be in a superposition of the states $ \ket{0} $ and $ \ket{1} $, i.e., $ \ket{q}  = a \ket{0} + b \ket{1} $ with $ a,b \in \CC $.
Hence, a qubit can be considered a vector in a two-dimensional complex vector space with the set $ \{  \ket{0},\ket{1}  \} $ forming a basis. 
A quantum system can be transformed from one state to another by the so called unitary transformations, which are reversible (see Fig.~\ref{Fig_quantum_gates_single_bit}, for some examples).  Once a classical measurement is taken the state of the qubit $ \ket{q}$ is changed from the superposition of  the states $ \ket{0} $ and $ \ket{1} $ to either $ 0 $ or $ 1 $, with respective probabilities of $ |a|^2 $ and $ |b|^2$. This measurement step is a non-unitary and non-reversible transformation.

We give a few examples to illustrate the working of a single qubit system. 
So, mathematically, if the state of a single qubit  $ \ket{q}  = a \ket{0} + b \ket{1} $ is represented by the matrix 
$
\left[\begin{smallmatrix}
a\\ 
b
\end{smallmatrix}\right]  $ and if a unitary transformation is represented by the $ 2 \times 2 $  unitary matrix  $ U = \mat{\alpha}{\beta}{\gamma}{\delta} $ then the resulting state is given by  $ Uq = \mat{\alpha}{\beta}{\gamma}{\delta} \left[\begin{smallmatrix}
a\\ 
b
\end{smallmatrix}\right]  $. 
Suppose the unitary matrix is $ X = \mat{0}{1}{1}{0} $ then the result of action of the unitary transformation $ X $ is  $ \mat{0}{1}{1}{0} \left[\begin{smallmatrix}
a\\ 
b
\end{smallmatrix}\right] =  \left[\begin{smallmatrix}
b\\ 
a
\end{smallmatrix}\right]  $. 
The unitary matrix $ X $ represents the quantum \textit{NOT} gate as it swaps the computational basis, i.e., it sends $ \ket{0} \to \ket{1} $ and $ \ket{1} \to \ket{0} $, see Fig.~\ref{Fig_quantum_gates_single_bit}.
Another very useful quantum gate is the Hadamard gate, whose matrix is given by $ H = \frac{1}{\sqrt{2}}\mat{1}{1}{1}{-1} $. It is often used to condition the input qubit which is often $ \ket{0} $ to a uniform superimposed state, $ H \ket{0} = \frac{1}{\sqrt{2}} (\ket{0} + \ket{1}) $. The above discussion can be easily extended to the multiple qubit systems. For example, a two qubit system has $ 4 $  possible states represented by $ \ket{00}, \ket{01},\ket{10} $ and $ \ket{11} $.

\begin{figure}[t]
\centering
\hspace{-2cm}
	\Qcircuit @C=1em @R=.7em {
		\lstick{a \ket{0} + b \ket{1}} 	& \qw & \gate{X} & \qw & \qw & \rstick{ b \ket{0} + a \ket{1}} \\
		\lstick{a \ket{0} + b \ket{1}} 	& \qw & \gate{H} & \qw & \qw & \rstick{  \frac{a}{\sqrt{2}}(\ket{0} +  \ket{1}) + \frac{b}{\sqrt{2}}(\ket{0} -  \ket{1})  }
	}

\caption {Quantum NOT and Hadamard gates} \label{Fig_quantum_gates_single_bit}
\end{figure}

\subsubsection{\textbf{Quantum search: Grover's algorithm }}
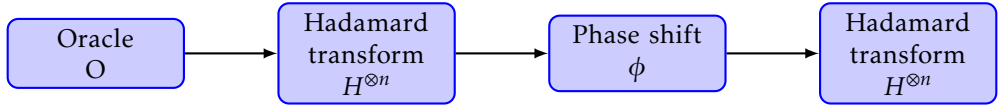
\begin{figure}
\centering
	\begin{tikzpicture}[%
	scale= 0.8,
	block/.style={
		rectangle,
		draw=blue,
		thick,
		fill=blue!20,
		text width=6em,
		align=center,
		rounded corners,
		minimum height=2em
	},
	block1/.style={
		rectangle,
		draw=blue,
		thick,
		fill=blue!20,
		text width=5em,
		align=center,
		rounded corners,
		minimum height=2em
	},
	line/.style={
		draw,thick,
		-latex',
		shorten >=2pt
	},
	cloud/.style={
		draw=red,
		thick,
		ellipse,
		fill=red!20,
		minimum height=1em
	}
	]
	\path (-9,0) node[block] (G) {Oracle \\ O}
	(-4.5,0) node[block] (H) {Hadamard transform \\ $ H^{\otimes n} $ }   
	(0,0) node[block] (I) {Phase shift \\ $\phi$}  
	(4.5,0) node[block] (J) {Hadamard transform \\ $ H^{\otimes n} $};    
	\draw[-latex,thick] (G) -- (H);                                   
	\draw[-latex,thick] (H) -- (I);  
	\draw[-latex,thick]   (I)--(J); 
	\end{tikzpicture}
	\caption{The Grover operator $ G $} \label{fig_Grover's_algo}
\end{figure}

We will briefly describe Grover's quantum search algorithm \cite{grover1996fast} for the sake of clarity and continuity of presentation in this work. 
Suppose that there is a list of $ N $ items, and we are required to search for one or more marked items. Also assume that we have access to an oracle, who when presented with an item from the list will answer if the item is marked or not. For example, the oracle can be assumed to be a black box function such that for any item with index $ x $ in the list, i.e., with $ 0 \leq x \leq N-1 $, we have  $ f(x)=1 $ if $ x $ is a marked item, otherwise $ f(x)=0 $. It is clear that with the classical computation total $ O(N) $ calls to oracle will be required for the search problem. However, with the quantum computation only $ O(\sqrt{N}) $ oracle calls are needed to solve the search problem with probability $ 1 $. Furthermore, if there are $ M $  solutions of the search problem (of size $ N $) then only $ O(\sqrt{N/M}) $  calls to the oracle is sufficient to solve the search problem.

\begin{figure}
\centering
	\begin{tikzpicture}[line cap=round,line join=round,>=triangle 45,x=0.5535055350553507cm,y=0.9385665529010234cm]
	\def\shift{3} 
	\def\down{-4} 
	\def \x {7}
	\draw [-] (0.,0.) -- (0.,2.);
	\draw [->] (0.,0.) -- (0.,1.);
	\draw [->] (1.,0.) -- (1.,1.);
	\draw [->] (2.,0.) -- (2.,1.);
	\draw [->] (3.,0.) -- (3.,1.);
	\draw [->] (4.,0.) -- (4.,1.);
	\draw [->] (5.,0.) -- (5.,1.);
	\draw [->] (6.,0.) -- (6.,1.);
	\draw [->] (7.,0.) -- (7.,1.);
	\draw (0.,0.)-- (9.,0.);
	\draw (-1,1.2) node[anchor=north] {$ \frac{1}{\sqrt{N}} $};
	\draw (0,-0.12) node[anchor=north] {0};
	\draw (1,-0.12) node[anchor=north] {1};
	\draw (6.75,-0.12) node[anchor=north west] {$N-1$};
	\draw (3.76,-0.12) node[anchor=north west] {$k$};
	\begin{scriptsize}
	\end{scriptsize}
	\draw [-] (\shift+ 10.,0.) -- (\shift+10.,2.);
	\draw [->] (\shift+10.,0.) -- (\shift+10.,1.);
	\draw [->] (\shift+11.,0.) -- (\shift+11.,1.);
	\draw [->] (\shift+12.,0.) -- (\shift+12.,1.);
	\draw [->] (\shift+13.,0.) -- (\shift+13.,1.);
	\draw [->] (\shift+14.,0.) -- (\shift+14.,-1.);
	\draw [->] (\shift+15.,0.) -- (1\shift+5.,1.);
	\draw [->] (\shift+16.,0.) -- (\shift+16.,1.);
	\draw [->] (\shift+17.,0.) -- (\shift+17.,1.);
	\draw (\shift+10.,0.)-- (\shift+19.,0.);
	\draw (\shift+10,-0.12) node[anchor=north] {0};
	\draw (\shift+11,-0.12) node[anchor=north] {1};
	\draw (\shift+16.75,-0.12) node[anchor=north west] {$N-1$};
	\draw (\shift+17.0,1.3) node[anchor=north west] {mean};
	\draw [dotted] (\shift+10.,0.9) -- (\shift+19.,0.9);
	\draw (\shift+13.9,-0.12) node[anchor=north west] {$k$};
	\draw [-] (\x + 0.,\down + 0.) -- (\x + 0.,\down + 2.);
	\draw [->] (\x + 0.,\down + 0.) -- (\x + 0.,\down + 1.);
	\draw [->] (\x + 1.,\down + 0.) -- (\x + 1.,\down + 1.);
	\draw [->] (\x + 2.,\down + 0.) -- (\x + 2.,\down + 1.);
	\draw [->] (\x + 3.,\down + 0.) -- (\x + 3.,\down + 1.);
	\draw [->] (\x + 4.,\down + 0.) -- (\x + 4.,\down + 3.0);
	\draw [->] (\x + 5.,\down + 0.) -- (\x + 5.,\down + 1.);
	\draw [->] (\x + 6.,\down + 0.) -- (\x + 6.,\down + 1.);
	\draw [->] (\x + 7.,\down + 0.) -- (\x + 7.,\down + 1.);
	\draw (\x + 0.,\down + 0.)-- (\x + 9.,\down + 0.);
	\draw (\x + -1,\down + 1.2) node[anchor=north] {$ \frac{1}{\sqrt{N}} $};
	\draw (\x + 0,\down -0.12) node[anchor=north] {0};
	\draw (\x + 1,\down -0.12) node[anchor=north] {1};
	\draw (\x + 6.75,\down -0.12) node[anchor=north west] {$N-1$};
	\draw (\x + 3.76,\down -0.12) node[anchor=north west] {$k$};
	\draw (\x + 4.1,\down + 3) node[anchor=north west] {$\frac{3}{\sqrt{N}}$};
	\end{tikzpicture}
	\caption {Amplitude amplification in one iteration of the application of Grover operator. Note that initially all the states have the same amplitude, but at the end the amplitude of the marked item is amplified.} \label{Fig_Grover_operator_amplitude}
\end{figure}
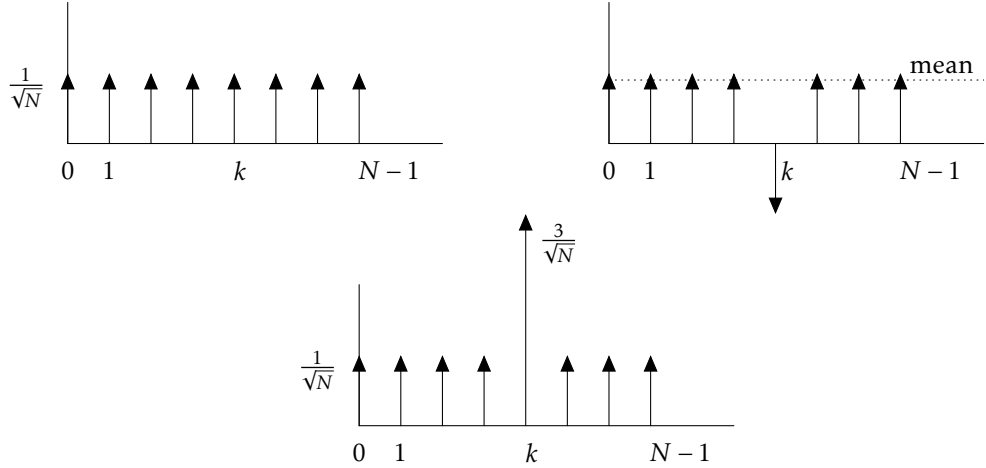

For simplicity, let us assume that $ N = 2^n $. Therefore, the discrete search space can be represented by a $ n $ qubit system. The algorithm starts with the system in $ \ket{0}^{\otimes n} $ state, i.e., all the $ n $ qubits are initially in the state $ \ket{0} $.  The action of the oracle $ O $ (see Eq.~$\! 6.3 $ in \cite{nielsen2002quantum} ) on the input $ \ket{x} $ can be represented by  
\begin{align}
O \ket{x} = (-1)^{f(x)} \ket{x}.
\end{align} 
The action of the Hadamard transform $ H^{\otimes n} $ on the input  $ \ket{0}^{\otimes n} $  is used to transform the input to a uniform superimposed state, say $ \ket{\psi} $, such that  
\begin{align}
\ket{\psi} = H^{\otimes n} (\ket{0}^{\otimes n} ) = \frac{1}{\sqrt{N}} \sum_{x=0}^{N-1} \ket{x}.
\end{align}
Next the algorithm proceeds with the repeated application of Grover operator $ G $ which is shown in the \mfig{fig_Grover's_algo}.  The third step in the \mfig{fig_Grover's_algo} is the phase shift which is transforming all non-zero basis state $ \ket{x} $ to $ -\ket{x} $ with $ \ket{0} $ remaining unchanged. It is easy to see that the phase shift $ \phi $ is the unitary operator $ (2 \ket{0} \bra{0} - I )$. It follows that the Grover operator $ G $ is essentially the unitary transformation 
\begin{align} \label{Eq_Grover_operator}
G = \left( H^{\otimes n} (2 \ket{0} \bra{0} - I) H^{\otimes n} \right) O = \left(2 \ket{\psi} \bra{\psi} - I  \right) O.
\end{align}
We note that at the heart of Grover's algorithm lies the operator $ G $. In fact, as noted earlier  $ G $ consists of $ \left(2 \ket{\psi} \bra{\psi} - I  \right) O$. The action of oracle amounts to changing the phase of the marked item. On the other hand it can be checked that 
\begin{align} \left(2 \ket{\psi} \bra{\psi} - I  \right) (\sum_{k=0}^{N-1} \alpha_k  \ket{k})  = \sum_{k=0}^{N-1} \left( -\alpha_k + 2 m \right) \ket{k} \end{align}
  where $ m = \frac{1}{N} \sum_{k=0}^{N-1} \alpha_k $ is the mean of $ \alpha_k $.  Hence, it is clear that the action of $ \left(2 \ket{\psi} \bra{\psi} - I  \right) $ is essentially an inversion about the mean. Therefore, one iteration of Grover operator results in amplification of the amplitude of the marked item (See \mfig{Fig_Grover_operator_amplitude}). 

\definecolor{qqqqff}{rgb}{0.,0.,1.}
\definecolor{xdxdff}{rgb}{0.49019607843137253,0.49019607843137253,1.}
\definecolor{uuuuuu}{rgb}{0.26666666666666666,0.26666666666666666,0.26666666666666666}

We note that if $ M=1 $ then the Grover operator $ G $ is applied $ R \approx \frac{\pi}{4} \sqrt{N} $ times. After $ R $ iteration of $ G $ the measurement of the system yields the marked item. To summarize Grover's algorithm for the case $ M=1 $ is as follows.
\begin{enumerate}
	\item Initialize the input $ n $ qubits as $ \ket{0}^{\otimes n} $.
	\item Apply $ H^{\otimes n} $
	\item Apply Grover operator $ G $ a total of $R \approx \frac{\pi}{4} \sqrt{N}$ times
	\item Measure the $ n $ qubits. \\
\end{enumerate}

\noindent \textbf{Some remarks:} 
\begin{enumerate}
	\item It is important to note that there is an optimum number of times (the integer $ R $ as above) that the Grover's iteration should be carried out to obtain the maximum probability of success (close to $ 1 $). This is in contrast to the classical behavior where the more iterations lead to better results.
	\item The Grover operator $ G $ (see \meqref{Eq_Grover_operator}) could be very nicely explained as rotation in a two dimensional space spanned by the vectors $ \ket{\alpha} $  and $ \ket{\beta} $ with $$ \ket{\alpha} = \frac{1}{\sqrt{N-M}} \sum_{x,\, x \text{ not marked}} \ket{x},\,\text{and }\,
	\ket{\beta} = \frac{1}{\sqrt{M}} \sum_{x,\, x \text{ marked}} \ket{x}. $$  
	The initial state $ \ket{\psi}  $  is in the span of $ \ket{\alpha} $ and $ \ket{\beta} $, as it is easily verified that $ \ket{\psi}  = \sqrt{\frac{N-M}{N}} \ket{\alpha} + \sqrt{\frac{M}{N}} \ket{\beta} $. If we set $ \cos(\frac{\theta}{2}) = \sqrt{\frac{N-M}{N}}$ then the  Grover operator may be written as 
	$
	\mat{\cos(\theta)}{-\sin(\theta)}{\sin(\theta)}{\cos(\theta)}.
	$
	The number of iteration $ R $ is chosen such that $ R \theta  $  is the angle required to turn the state vector $ \psi $ very close to $ \ket{\beta} $. For example, if $ 1= M <<N $ then we get $ \theta \approx \sin(\theta) = 2 \sin(\frac{\theta}{2}) \cos(\frac{\theta}{2}) \approx 2 \sqrt{\frac{1}{N}} $. Also as $ \psi $ is initially almost along the $ \ket{\alpha} $ the total turn is approximately $ \frac{\pi}{2} $. Therefore, we must have $ R \theta \approx 2 R \sqrt{\frac{1}{N}}  \approx \frac{\pi}{2} $. And we see that  $R \approx \frac{\pi}{4} \sqrt{N}$.  We refer readers to  \cite{nielsen2002quantum}  for further details.
\end{enumerate}

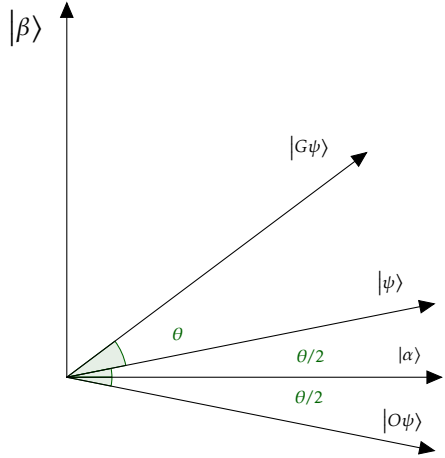
\begin{figure}
	\centering
	\definecolor{qqwuqq}{rgb}{0.,0.39215686274509803,0.}
	\begin{tikzpicture}[line cap=round,line join=round,>=triangle 45,x=1.0cm,y=1.0cm]
	\draw [shift={(0.,0.)},color=qqwuqq,fill=qqwuqq,fill opacity=0.10000000149011612] (0,0) -- (11.309932474020219:0.8) arc (11.309932474020219:36.86989764584403:0.8) -- cycle;
	\draw [shift={(0.,0.)},color=qqwuqq,fill=qqwuqq,fill opacity=0.10000000149011612] (0,0) -- (-11.309932474020219:0.6) arc (-11.309932474020219:0.:0.6) -- cycle;
	\draw [shift={(0.,0.)},color=qqwuqq,fill=qqwuqq,fill opacity=0.10000000149011612] (0,0) -- (0.:0.6) arc (0.:11.309932474020215:0.6) -- cycle;
	\draw [->] (0.,0.) -- (5.,0.);
	\draw [->] (0.,0.) -- (0.,5.);
	\draw [->] (0.,0.) -- (4.9029033784546,0.9805806756909203);
	\draw [->] (0.,0.) -- (4.9029033784546,-0.9805806756909203);
	\draw [->] (0.,0.) -- (4.,3.);
	\draw (-0.9,5) node[anchor=north west] {$\left| \beta\right\rangle$};
	\begin{scriptsize}
	\draw[color=black] (4.56,0.29) node {$\left| \alpha \right\rangle$};
	\draw[color=black] (4.3,1.25) node {$\left| \psi \right\rangle$};
	\draw[color=black] (4.46,-0.59) node {$\left| O \psi \right\rangle$};
	\draw[color=black] (3.22,3.05) node {$\left| G\psi \right\rangle$};
	\draw[color=qqwuqq] (1.48,0.57) node {$\theta$};
	\draw[color=qqwuqq] (3.22,-0.27) node {{\scriptsize ${\theta}/{2}$}};
	\draw[color=qqwuqq] (3.24,0.27) node {{\scriptsize ${\theta}/{2}$}};
	\end{scriptsize}
	\end{tikzpicture}
	\caption {Grover operator as rotation of the state vector $ \ket{\psi} $ towards the superposition $\ket{\beta} $ of all the solution vectors of the search problem.}
\end{figure}

\subsection{\textbf{Method IV: Quantum exhaustive search algorithm}} \label{subsect5_quantum_exhaustive}
Several quantum search algorithms based on Grover's algorithm have been treated in the literature. If one is looking for $ M $ marked things out of $ N $ things, with number $ M $  not known a priori, then the expected running time of these algorithms are $ O(\sqrt{\frac{N}{M}}) $. D\"urr and H\o yer  have given a quantum algorithm in \cite{durr1996quantum} to find the minimum cost. Algorithm~$ \ref{Algo_2} $ described below is an adaptation of their algorithm suited for our framework to solve a trajectory optimization problem.  We will give a computational example in \mref{subsect_exhaustive_quantum}. We note that Algorithm~$ \ref{Algo_2} $ finds the global minimum for the problem described in \mref{Discretization} with the cost function being the discrete version of  \meqref{Eq_Main_Cost}. It takes a list of size $ N $ as input. This list contains all the possible states of the system. The algorithm proceeds by first selecting a random state from the input list and setting its associated cost as the running minimum. The value of parameter $ \lambda $ and the variable $ m $ is also initialized. The variable $ m $ is initially set to $ 1$  and later in the algorithm it is adaptively scaled by $ \lambda $, by setting its value to $ \lambda m $ if the search for a state with lower associated cost than the current running minimum cost is not successful. The variable $ m $ effectively controls the number of  Grover's rotation $ r $ that is performed. A priori, it is not known that how many of the states have their associated costs less than the current running minimum cost and therefore $ m $ is adaptively scaled to ensure that a correct value for $ r $ is selected as required by Grover's search algorithm. 

\begin{algorithm}[ht]
	\SetKwInOut{Input}{Input}
	\SetKwInOut{Output}{Output}
\tcc{	Description: The objective of this quantum algorithm is find the minimum cost for the problem described in \mref{Discretization} with the cost function being the discrete version of  \meqref{Eq_Main_Cost}  and $ N = ST $ with $ S $ and $ T $ as in  \meqref{Eq_S}   and  \meqref{Eq_T}.}
	\underline{function findMinCostQuantum} ($ \{Y_0,Y_1,Y_2 \ldots ,Y_ {N-1}   \} $)\;
	\Input{a list of size $ N $ containing all the possible states of the system.}
	\Output{(the index of an admissible state which has the minimum cost, the associated minimum cost).}
	Set $ y = $  a random state from the input list.\\
	Set $ m=1 $. Set the total number of Grover's rotation, $ G = 0 $. \\
	Set the value of the parameter $\lambda  $ ($ \lambda = \frac{8}{7} $ is used in \cite{boyer1998tight}, but $ \lambda = 1.34 $ is better as suggested in \cite{bulger2003implementing}). \\
	\While { $ G $ is less than $ 22.5 \sqrt{N}  + 1.4 \log^2(N)$  }{ 
	Choose $ r$ uniformly at random from $ \{0,1,2, \ldots \ceil{m-1}  \} $ \\ 
	Perform Grover's search with $ r$ rotation to find a state $ x $ with $ cost(x) < cost(y) $. \\
	Increment the total number of Grover's rotation, $ G = G + r $.\\
	\qquad \uIf{Search is successful}
	{
		\qquad Set $ y = x $, $ m=1 $  \\
			}
	\qquad \uElse 
	{
		\qquad  Set $ x=y $, $ m = \lambda m $ \\
		}
	} 
     Return ($ x $, $ cost(x) $)
	\caption{Quantum algorithm for trajectory optimization} \label{Algo_2}
\end{algorithm}

\vspace{0.5cm}
The algorithm runs until the total accumulated number of Grover's rotation $ G $ is less than $ 22.5 \sqrt{N}  + 1.4 \log^2(N)$. We refer readers to \cite{durr1996quantum} for further details. 
We also note that the running time of the above algorithm is optimal up to a constant prefactor. See \cite{bulger2003implementing} and \cite{lara2014new} for other variants of the above algorithm with smaller constant prefactors. We note that in place of Algorithm~$ \ref{Algo_2}$, these other variants may also be employed to solve the trajectory optimization problem.

It is clear that the quantum algorithms discussed above are of the order of $ O(\sqrt{N}) $. Hence, they are far superior to their classical counterparts which are of the order of $ O(N) $.

\subsection{\textbf{Method V: Quantum random search algorithm}}\label{subsect5_quantum_random}
As before we consider the problem of finding the minimum cost for the problem described in \mref{Discretization}, with the cost function being the discrete version of  \meqref{Eq_Main_Cost}  and $ N = ST $ with $ S $ and $ T $ as in  \meqref{Eq_S}   and  \meqref{Eq_T}.
As the first step of this algorithm a predetermined number of states are selected uniformly at random out of the total $ N$ states. Then Algorithm~$ \ref{Algo_2} $ is applied on these selected states. Similar to the classical random search algorithm, this algorithm is especially useful in situations wherein the discrete search space is so big that an exhaustive search has a prohibitively high computational cost. We will consider a computational example in \mref{subsect_random_quantum}.

\subsection{\textbf{Method VI: Quantum hybrid search algorithm}} \label{subsect5_quantum_hybrid}
 We consider a quantum hybrid algorithm to solve the problem described in \mref{Discretization} with the cost function being the discrete version of  \meqref{Eq_Main_Cost}  and $ N = ST $ with $ S $ and $ T $ as in  \meqref{Eq_S}   and  \meqref{Eq_T}.
A quantum hybrid search algorithm is a combination of a quantum random search algorithm and a pure exhaustive quantum search algorithm. It combines the best attributes of both a quantum random search algorithm and a pure quantum exhaustive search algorithm. To begin with, in a hybrid algorithm  first a fixed number of states are selected uniformly at random out of the total $ N $ states and on these selected states the quantum algorithm described earlier (Algorithm~$ \ref{Algo_2} $) is applied. Next, using the approximate solution thus obtained a finer discretization is carried out. Finally, an exhaustive quantum search is carried out using Algorithm~$ \ref{Algo_2} $ in this finer grid. A computational example is given in \mref{subsect_hybrid_quantum}.

Here it is interesting to note that \cite{bulger2007combining} also describes a hybrid method, by effectively combining local search  with Grover's algorithm, to obtain the global optimum of an objective function. The quantum basin hopper algorithm, proposed in \cite{bulger2007combining}, considers a black-box real-valued objective function $ f $ defined on a discrete domain $ S $,  of size $ 2^n $. It also makes certain regularity assumptions on the objective function and its domain. As the method in \cite{bulger2007combining} does not deal with the problem of discretization, it is not directly applicable for solving a continuous trajectory optimization problem. It is in this context that our framework could be useful. Following our discretization approach, depending upon the presence of any regularity or local structure in the problem, various local search methods and variants of Grover's algorithms, including the method in \cite{bulger2007combining}, can be used in our framework. Although, in our present work (including the computational examples that follow), we do not make any assumptions about the availability of local search methods or about the presence of any local structure in the optimization problems.

\section{Computational examples}\label{Computational examples} 

Now we consider a couple of examples for illustrating our method. 
\subsection{\textbf{Brachistochrone problem}}
The first example that we take is the well-known Brachistochrone problem. The problem is to find the required trajectory of a particle starting from the rest under the influence of gravity, without any friction, such that it slides from the one fixed point to the other in the least possible time.
We let $ \vec{X}(t) = x(t) \vec{i} + y(t) \vec{j} $ to represent the position of the particle. The boundary condition that we consider is 
$ (x(\tau_0),y(\tau_0)) = (0,2) $ and $ (x(\tau_f),y(\tau_f)) = (\pi,0) $ with $ \tau_0 = 0 $. The goal is to minimize $ \tau_f $ given by
\begin{align}
\tau_f = \int_{0}^{\pi} \, \sqrt{\frac{1 + \left(\frac{dy}{dx}\right)^2 }{2gy}} \, dx.
\end{align}
Next we describe some possible approaches for discretization that could be employed to solve the brachistochrone problem.

\subsubsection{\textbf{Global dicretization in physical space}} \label{sub_sect_discrete_physical_space}

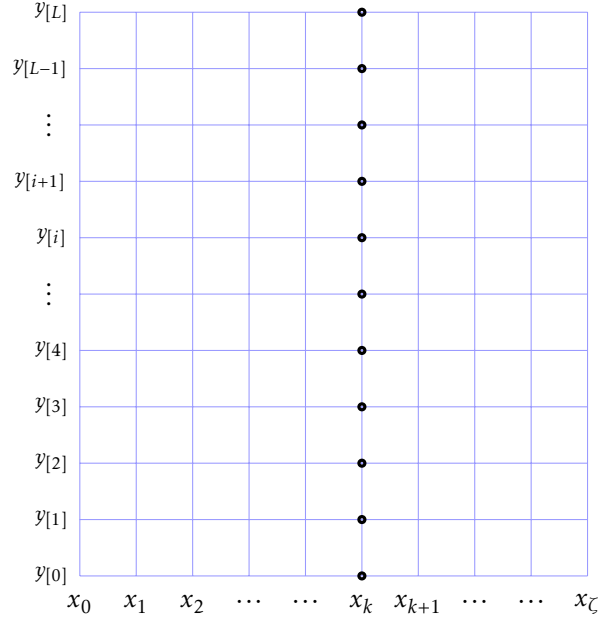
\begin{figure}
\centering	
\begin{tikzpicture} [scale=0.75]
\draw [step=1cm,blue!40,thin] (0,0) grid (9,10); 
\foreach \k in {0,1,2}
\node  at (\k,-0.5) {$ x_{\k} $};
\node  at (3,-0.5)  {$ \cdots  $};
\node  at (4,-0.5)  {$ \cdots  $};
\node  at (5,-0.5)  {$ x_k  $};
\node  at (6,-0.5)  {$ x_{k+1}  $};
\node  at (7,-0.5)  {$ \cdots  $};
\node  at (8,-0.5)  {$ \cdots  $};
\node  at (9,-0.5)  {$ x_\myN  $};
\foreach \k in {1,2,3,4}{
\pgfmathtruncatemacro \y {2*\k};
\node  at (-0.5,\k) {{\small $ y_{[\k]}$}};}
\node [rotate=90] at (-0.5,5)  {$ \cdots  $};
\node  at (-0.5,6)  {{\small $y_{[i]}    $}};
\node  at (-0.7,7)  {{\small $ y_{[i+1]}  $}};
\node [rotate=90] at (-0.5,8)  {$ \cdots  $};
\node at (-0.7,9) {{\small $ y_{[L-1]}  $}};
\node at (-0.5,10) {{\small $ y_{[L]}$}};
\node at (-0.5,0) {{\small $ y_{[0]} $}};
\foreach \i in {0,...,10}{
\draw [color=black, very thick] (5,\i) circle (1.5pt) ;
	}
\end{tikzpicture}
\caption{Brachistochrone problem: discretization of both the dependent variable $ y $ as well as the independent variable $ x $ is carried out. Note that at $ x=x_k $, the corresponding $ y(x_k) $ is allowed to take values only from the set  $ \{ y_{[i]} :\: i = 0,1 \cdots L \} $. The possible values for $ y(x_k) $ are shown with thick black dots.} \label{Fig_grid}
\end{figure}

Out of the infinitely many possible paths the goal is to pick the path with the minimum time. For practical purposes it is sufficient to find a `good enough' approximate solution (based on acceptable levels of errors).
The discretization is chosen based on what constitute  a `good enough' solution in a given context. In our present example, the physical space consists of the rectangle $ [0,\pi] \times [0,2] $ in $ \RR^2 $. Of course, one may as well consider a bigger rectangle to allow for the possibility of better solutions. One can now discretize the rectangle $ [0,\pi] \times [0,2] $ in several possible ways. For example, let $ x_k = \frac{k\pi}{\myN} $ for $k = 0,1 \cdots \myN $. Let $ y_{k} \in \{ y_{[i]} = \frac{2i}{L} :\: i = 0,1 \cdots L \} $ for $ k=1 $ to $ \myN-1 $. We set $ y_0 =2 $  and $ y_\myN = 0 $ to take into account the boundary condition. Next we use Lagrange interpolation to set $ y(x) = \sum_{k=0}^{\myN} y_{k} \phi_k(x) $. Here $ \phi_k(x)  $  is chosen such that 
$ \phi_k(x_j) = 1 $ if $ j=k $ and $ \phi_k(x_j) = 0$ otherwise. More explicitly, 
\begin{align}
\phi_k(x) = \prod_{j=0,~ j\neq k}^{\myN} \, \frac{(x-x_j)}{x_k-x_j}.  
\end{align}
Essentially, it means that once a discretization of $ x $ is carried out then we require that at $x= x_k $, the corresponding $ y(x_k) $ can only have values from the set $ \{y_k\} $ (see \mfig{Fig_grid}). We note that with the above discretization there are $ N = (L+1)^{\myN-1} $  total possible cases need to be considered.

Theoretically, one can always approach the correct solution to the desired accuracy by making the discrete grid finer. However, the discrete search space grows very rapidly and other methods should be combined with exhaustive search in practical applications. We consider a few examples to compare the correctness of the results with respect to the size of the discrete search spaces and the points chosen for discretization. Examples A and B are chosen to highlight the effect of uniform and non-uniform grid spacing respectively (for $ \myN=5 $). Example C is chosen to highlight the effect of having a finer grid resolution.

\begin{itemize}
\item \textbf{Example A:} Consider the path obtained by using the Lagrange interpolation on the $ (x,y) $ coordinates listed in Table \ref{Table_Example_A}, 
\begin{align*}
	f_A (x) &=  \sum_{i=0}^{\myN} a_i x^i		
	\end{align*}
with 
\begin{align*}
 a_0 &= 2,  \qquad & a_1 &= - 2.72552840044871, \\
 a_2 &=2.37999238659866, \qquad & a_3 &= - 1.34381393471665, \\
 a_4 &=0.387647767429487, \qquad & a_5 &=-0.0425490057689243.
\end{align*}
This path can be discovered by taking   $ L=40 $ and $ \myN=5 $ in the preceding discussion. Therefore, the total size of the discrete search space is $ 41^4 = 2825761 $. 
We note that the time obtained in this case, $ 1.00946330885 $ sec, differs from the correct analytical solutions $ 1.0035449615773016 $ sec by only about  $0.6$\%. \\

\item \textbf{Example B:} Let us now consider a different discretization and pick the points as shown in Table \ref{Table_Example_B}. The path obtained by using the Lagrange interpolation in this case is  
\begin{align*}
f_B (x) &=  \sum_{i=0}^{\myN} a_i x^i	,	
\end{align*}
with 
\begin{align*}
 a_0 &= 2,  \qquad & a_1 &= - 2.72732233702330, \\
 a_2 &=2.08090516797775, \qquad & a_3 &= - 0.962129500276245, \\
 a_4 &=0.229887380333958, \qquad & a_5 &=-0.0213407110121635.
\end{align*}
This path can be discovered by using a discrete search space of size $ 41^4 = 2825761 $. The minimum time obtained is  $ 1.00960568357 $ sec. The error in the time obtained in this case is $ 0.6 $\%, the same as the previous case.

\item \textbf{Example C:} In this example we get a more accurate result by considering more points as shown in Table \ref{Table_Example_C}. The size of the discrete search space to discover this path, if we search exhaustively, is $ 201^{12} $. The minimum time obtained in this case is $ 1.00436810144 $ sec and the error percentage is about $ 0.08 $\%.
\end{itemize}

\begin{table}
	\begin{center}
		\small\setlength{\tabcolsep}{5pt}
		\begin{tabular}{cccc}
			\toprule \toprule
			& \multicolumn{3}{c}{Root Mean Square and Percentage Error} \\ 
			\hline
			\midrule
			 & Example A & Example B & Example C \\
			\midrule 
		RMSE:  & $ 0.0431602981808 $ & $  0.0441659854713 $ &  $ 0.00839381927831 $  \\
		\% Error: & $  0.6 $ & $ 0.6 $ & $  0.08  $ \\
			\bottomrule
		\end{tabular} 
		\caption{Root mean square and the percentage errors between the minimum time obtained in Example A, B and C and the correct analytical minimum time for the brachistochrone problem.} \label{Table_Example_RMSE} 
	\end{center}
	\end {table}

\begin{figure}
	\centering
	\includegraphics[scale=0.75]{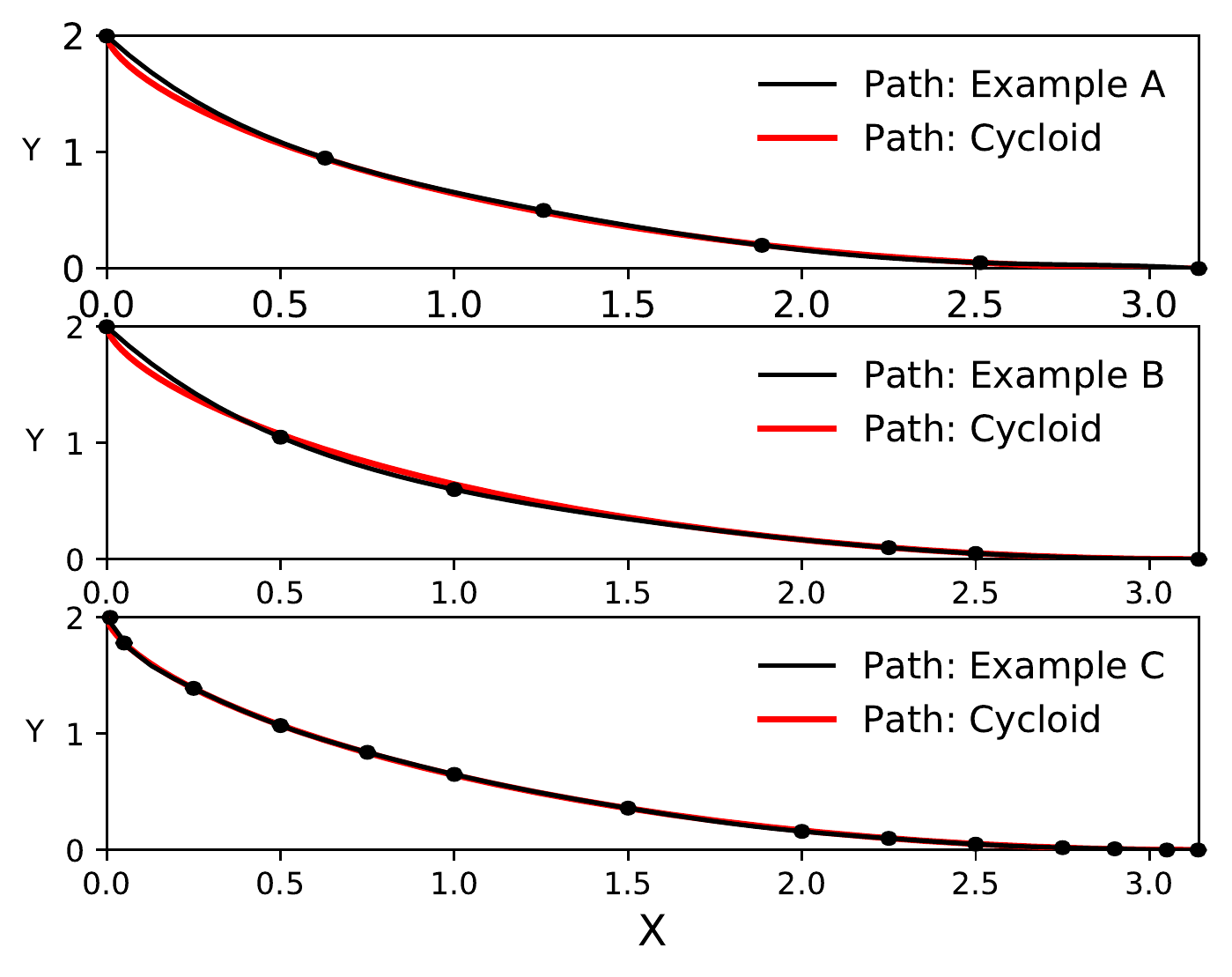}
	\caption{The paths for Examples A, B and C are plotted along with the cycloid resulting from the correct analytical solution of the brachistochrone problem. The nodes in the above examples (Examples A \& B: 6 nodes, Example C: 14 nodes) represent the solution points on the discretization grid (see Table \ref{Table_Example_A}, Table \ref{Table_Example_B}, and Table \ref{Table_Example_C}).}\label{Fig_examples}
\end{figure}

In \mfig{Fig_examples} the paths for the above examples are plotted along with the cycloid resulting from the correct analytical solution. Errors associated with these paths are shown in Table~\ref{Table_Example_RMSE}. We note that results shown in \mfig{Fig_examples} for Examples A, B and C are in good agreement with the analytical solution. As expected, Example C has the least error (as noted in Table~\ref{Table_Example_RMSE}) which can be attributed to a finer grid resolution. 

\begin{figure}[ht]
	\centering
 \label{Fig_errors}
	\includegraphics[scale=0.75]{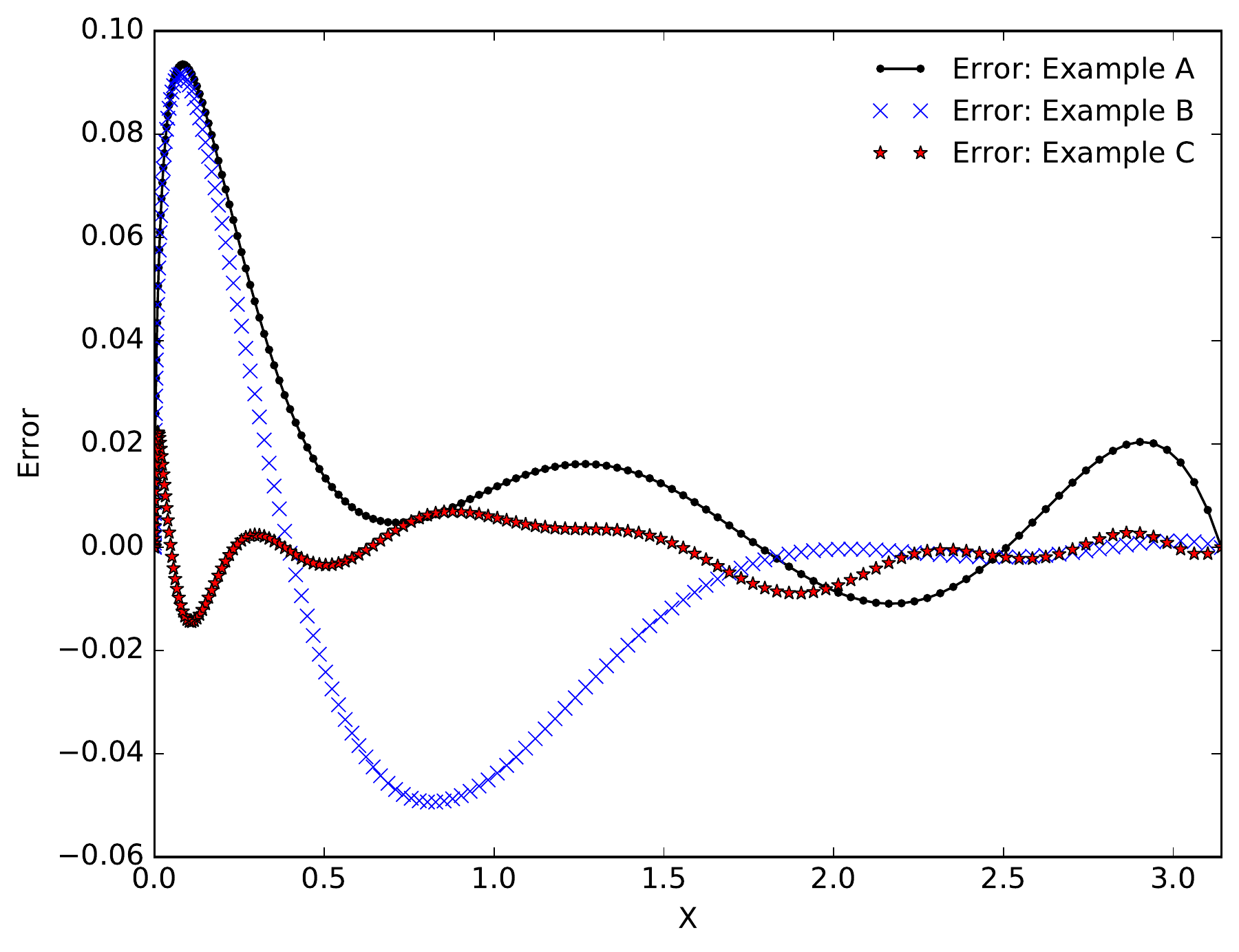}
	\caption{ Errors in the paths obtained in the Example A, B and C and the cycloid resulting from the correct analytical solution of the brachistochrone problem.}
\end{figure}

\begin{table}
	\begin{center}
		\small\setlength{\tabcolsep}{5pt}
		\begin{tabular}{c*{6}{c}}
			\toprule \toprule
			& \multicolumn{6}{c}{Example A} \\ 
			\hline
			\midrule
			x:& $ 0 $  & $ \dfrac{\pi}{5} $ &  $ \dfrac{2\pi}{5}$ & $ \dfrac{3\pi}{5} $ & $ \dfrac{4\pi}{5} $  & $ \pi $ \\
			\midrule 
			y: & $ 2 $ & $ 0.95 $ &  $ 0.50 $ &  $ 0.20 $ &  $ 0.05 $ & $ 0 $ \\
			\bottomrule
		\end{tabular} 
		\caption{Uniformly spaced sample points in $ x $ for the brachistochrone problem.} \label{Table_Example_A}
	\end{center}
	\end {table}
	
 	\begin{table}
 		\begin{center}
 			\small\setlength{\tabcolsep}{5pt}
 			\begin{tabular}{c*{6}{c}}
 				\toprule \toprule
 				& \multicolumn{6}{c}{Example B} \\ 
 				\hline
 				\midrule
 				x:& $ 0 $  & $ 0.5 $ & $ 1.0 $ & $ 2.25 $ & $ 2.5 $ & $ \pi $ \\
 				\midrule 
 				y: & $ 2 $ & $ 1.05 $ &  $ 0.60 $ &  $ 0.10 $ &  $ 0.05 $ & $ 0 $ \\
 				\bottomrule
 			\end{tabular} 
 			\caption{Non-uniformly spaced sample points in $ x $ for the brachistochrone problem.} \label{Table_Example_B}
 		\end{center}
 		\end {table}

\begin{table}
	\begin{center}
		\small\setlength{\tabcolsep}{3pt}
		\begin{tabular}{c*{14}{c}}
			\toprule \toprule
			& \multicolumn{14}{c}{Example C} \\ 
			\hline
			\midrule
			x:& $ 0 $ & $ 0.05 $ & $ 0.25 $ & $ 0.5 $ & $ 0.75 $ & $ 1.0 $ & $ 1.5 $ &  $ 2.0 $ & $ 2.25 $ & $ 2.5 $ & $ 2.75 $ & $ 2.9 $ & $ 3.05 $ & $ 3.14 $ \\
			\midrule 
			y: & $ 2 $ & $ 1.78 $ & $ 1.39 $ & $ 1.07 $ & $ 0.84 $ & $ 0.65 $ & $ 0.36 $ & $ 0.16 $ & $ 0.1 $ & $ 0.05 $ & $ 0.02 $ & $ 0.01 $ & $ 0 $ & $ 0 $ \\
			\bottomrule
		\end{tabular} 
		\caption{A finer grid resolution for obtaining a more accurate solution for the brachistochrone problem.} \label{Table_Example_C}
	\end{center}
	\end {table}

\noindent \textbf{Some remarks:} 
\begin{enumerate}
	\item One can easily put the above description in the general discretization framework that we discussed in the \mref{Discretization}. For example,  setting $ x(t) =t $ and $ y(t) = f(t) $ is the same as saying $ y = f(x) $. 
	\item In the above discussion, we have picked the discrete search spaces rather naively, just to show that it is possible to get a good solution (within $ 0.6 \% $ accuracy) even in such cases. One can get a smaller discrete search space by starting with an initial guess, for example in brachistochrone problem this could be a linear path, or one can do multiple passes of search by progressively refining the discrete grid in later passes.  
\end{enumerate}

\subsubsection{\textbf{Discretization in coefficient space}}
Sometimes due to the nature of the optimization problem it is more convenient and computationally efficient to discretize in the coefficient space. We will solve the brachistochrone problem by employing this technique.  Suppose  $ y(x) = \sum_{n=0}^{\myN-1} a_n x^n $ is a possible path. We want to discretize the coefficients $ a_n $. In order to perform this discretization, we need to first determine the range of the coefficients $ a_n $, i.e., their minimum and maximum values. The boundary condition, $ y(0) =2 $ implies  
\begin{align}\label{eq_boundary_a_0}
a_0 =2,
\end{align}
and we use the boundary condition $y(\pi) =0 $ to fix 
\begin{align} 
 a_1 = \frac{1}{\pi} \left(-a_0 - \sum_{n=2}^{\myN-1} -a_n\pi^n \right).
\end{align}
	Now guided by the physical space of this problem we impose the condition that 
	\begin{align} \label{eq_y_constraint}\tag{6.5}
	| y(x) | \,\, \leq \, 2  \qquad \text{for }\,  x \in [0,\pi].
	\end{align}
	
	Further, we note some interesting results in approximation theory (ref.~\cite{rahman2002analytic}, Theorems $ 16.3.1 $ and $ 16.3.2 $ therein) on the bounds on coefficients of bounded polynomials.  According to these theorems, given a bound on $ y(x) $, such as in \meqref{eq_y_constraint}, the bound on the coefficients of the polynomial $ y(x) $ can be determined by using the coefficients of Chebyshev polynomials of the first kind. In the following, we provide more details on these theorems to demonstrate their application in our context.
	
	Let $ T_n(x) $ be the Chebyshev polynomial of the first kind of degree $ n $ with the series expansion given as $ T_n(x) = \sum_{k=0}^{n}  t_{n,k} x^{k}  $. Let $ \frac{1}{2} y(\frac{\pi}{2} (x+1)) = \sum_{n=0}^{\myN-1} b_n x^n $. It is clear that $ |\frac{1}{2} y(\frac{\pi}{2} (x+1))| \leq 1 $ for  $ |x| \leq 1  $, and therefore it satisfies the hypothesis of the earlier mentioned theorems in \cite{rahman2002analytic}. Hence, we obtain the following bounds on the coefficients of $ \frac{1}{2} y(\frac{\pi}{2} (x+1)) $
	\begin{align} \label{eq_coefficient_bound1}
	|b_{n-2 \mu}|  &\leq |t_{n,n-2\mu}| \qquad \text{for $ \mu =0,\cdots, \lfloor{\frac{n}{2}}\rfloor $ }, \tag{6.6}\\
	|b_{n-2 \mu-1}| & \leq  |t_{n-1,n-2\mu-1}| \qquad \text{for $ \mu =0,\cdots, \lfloor{\frac{n-1}{2}}\rfloor $ }. \label{eq_coefficient_bound2} \tag{6.7}
	\end{align}  
	Of course, the coefficients $ a_1,\cdots a_{\myN-1} $ can be expressed in terms of the coefficients $ b_1,\cdots b_{\myN-1}  $ using the relation $ \sum_{n=0}^{\myN-1} b_n x^n  = \frac{1}{2} \sum_{n=0}^{\myN-1} a_n \left(\frac{\pi}{2}(x+1)\right)^n $, and therefore from the above relations using the bounds on coefficients $ b_1,\cdots b_{\myN-1}  $ coefficients $ a_1,\cdots a_{\myN-1} $ can be discretized.

\noindent \textbf{Example D:}
We set $ \myN=6 $ and note that $ T_5(x) = 16 x^5 -20x^3 +5x$ and $ T_4(x)=8x^4-8x^2+1 $. On using \meqref{eq_boundary_a_0} through \meqref{eq_coefficient_bound2}, we obtain the following conditions 
	\begin{alignat}{3}
	&|b_2| \leq 8, \qquad |b_3| \leq 20, \qquad |b_4| \leq 8, \qquad |b_5| \leq 16,  \nonumber\\
	&a_5 = 64b_5/\pi^5, \quad a_4 = 32 (b_4 - 5 b_5)/\pi^4, \quad a_3 = 16 (b_3 - 4b_4 + 10 b_5)/\pi^3   \nonumber \\ 
	&a_2 = 8(b_2 - 3b_3 + 6b_4 - 10b_5)/\pi^2,  \quad a_1 = -2/{\pi} - a_2 \pi -a_3 \pi^2 - a_4 \pi^3 - a_5 \pi^4, \quad a_0 = 2.\tag{6.8}
	\end{alignat}
	
	Using the above relations, we discretized the coefficients $ a_2,a_3,a_4 $ and $ a_5 $ by letting each of $ b_2 $ and $ b_4 $ to take values in the set  $ \{ -8,-7.8, \cdots 7.8,8\}$  of cardinality  $81$, $ b_3 $ to take values in the set $  \{-20,-19.5, \cdots,\allowbreak 19.5,20 \}$  of cardinality  $81$,  and $ b_5 $ to take values in the set $  \{-16,-15.8, \cdots, 15.8,16 \}$  of cardinality  $161$. Then one can find the  path $ y(x) =-0.0418 x^5 + 0.3942 x^4 - 1.4449 x^3 + 2.7559 x^2 - 3.1831 x + 2 $ with the corresponding minimum time of $ 1.01067460259 $ sec. We note that the error in the minimum time calculated is about $ 0.7\% $ and the total number of paths considered is $  81 \times 81 \times 81 \times 161=  85562001 $.  
	Of course, the accuracy can be further improved by making the discretization grid finer as well as increasing $ \myN $. \\
	
\noindent \textbf{A remark:} In the above example, an error of $ 0.7\% $ in the minimum time was obtained with a search space of the size $ 85562001 $. A much smaller search space will suffice, if a larger error in calculating the minimum time is allowed. For comparison we note that, using a similar approach as discussed above for a degree three polynomial, i.e., for $ \myN=4 $, we obtained the minimum time with an error of $ 1.4\% $, by considering a search space of the size $ 1600 $.

\subsubsection{\textbf{Classical random search algorithm in physical space}} \label{subsect_random}

\noindent \textbf{Example E:}
Now we consider the same global discretization using the actual physical space that we considered in the Example A, \mref{sub_sect_discrete_physical_space}, with $ \myN=5 $ and $ L=40 $.
Then the best path that we obtained on a trial run was
\begin{align*}
f_E (x) &=  \sum_{i=0}^{\myN} a_i x^i,		
\end{align*}
with 
\begin{align*}
 a_0 & = 2,    &a_1 &= -  3.09025847836764, \\
 a_2 &=2.99741834941916,  &a_3 &= -1.62937439584393, \\
 a_4 &=0.427749260611848,  &a_5 & =0.0425490057689243,
\end{align*}
 and the corresponding minimum time found was $ 1.00993422472 $ sec. We note that the error in the minimum time calculated is about $ 0.7\% $ and the total number of paths considered is $ 5000 $.  
Of course, the accuracy can be increased by making the discretization grid finer as well as increasing the number of paths considered.

\begin{table}
\caption {A comparison of performances of different search methods for the brachistochrone problem. (Assume $ \epsilon = 2.46 $.)} \label{tab:title} 
\begin{center}
\begin{tabular}{lllll}\toprule \toprule
{\small Search }& {\small 	Search method description} & {\small Cost} & {\small Minimum time}  & {\small  Error} \% \\ 
{\small method} & & & & \\
	\hline  \\
{\small I} &	Classical exhaustive & 2825761  & 1.0095 &  0.6 \\ 
&	{\small (Example A, physical space)} & & & \\ \\
{\small I} &	Classical exhaustive   &   85562001   &  1.0107 &  0.7 \\ 
&	{\small (Example D, coefficient space)} & & & \\ \\
{\small II} &	Classical randomized & 5000  & 1.0099  & 0.7 \\ \\
{\small III}&    Classical hybrid   	&  8249 &  1.0085 & 0.5 \\ \\
{\small IV}&	Quantum exhaustive & $  \sqrt{2825761} \, \epsilon \approxeq 4135$ & 1.0095 & 0.6 \\ \\ 
{\small V}&	 Quantum random & $  \sqrt{5000} \, \epsilon \approxeq 174$ & 1.0099 & 0.7 \\ \\ 
{\small VI}&	Quantum hybrid   & {\footnotesize $ (\sqrt{5000} + \sqrt{3249}) \, \epsilon  \approxeq 314 $} & 1.0085 & 0.5  \\ \\
	\bottomrule
\end{tabular} 
\end{center}
\end{table}

\begin{figure}[ht] 
	\centering
	\includegraphics[scale=0.7]{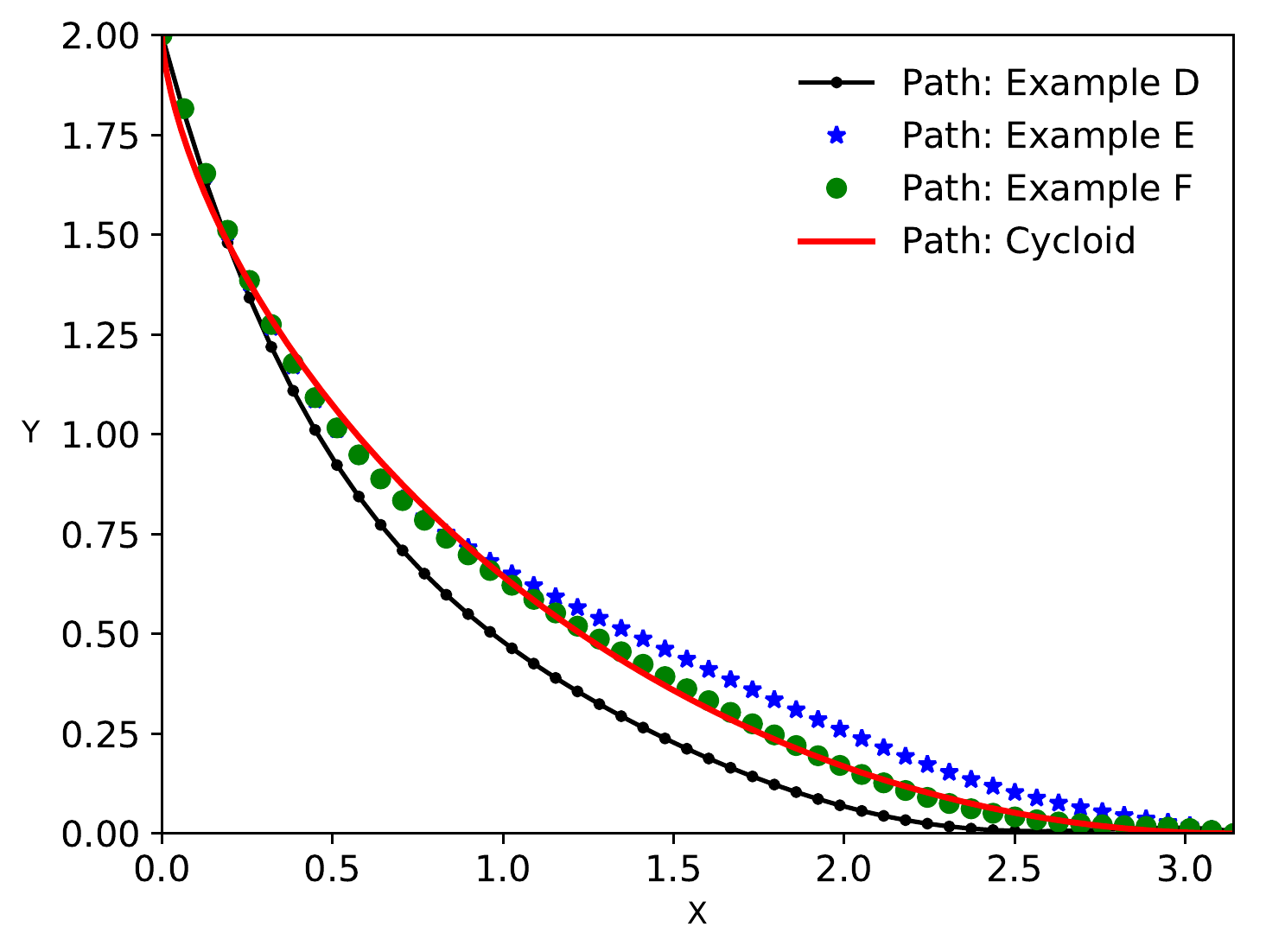}
	\caption{Comparison of optimal paths obtained by using different approaches such as a discretization in the coefficient spaces (Example D), random search in a discretized physical space (Example E) and a hybrid search combining global exhaustive search with coarse-grained random search (Example F).}\label{Fig_oeff_ran_ex}
\end{figure}

\subsubsection{\textbf{Classical hybrid search algorithm}} \label{subsect_hybrid}
In many situations, it is a good idea to follow the chosen discretization method by a random search to obtain an approximate solution. Then one can pick a finer discrete grid to further improve the result of the random search by performing an exhaustive search in a finer discrete grid near the approximate solution given by the random search.
 
\noindent \textbf{Example F:} The result obtained from the random search in the \mref{subsect_random} is further refined by performing an exhaustive search near the best path given by the random search algorithm. We consider  $ x_k = \frac{k\pi}{\myN} $  for $ k=0 $ to $ \myN $ and $ \myN=5 $. Further, we set $ y(x_0) = 2 $, $ y(x_1) \in \{0.8,0.85,0.90$, $0.95\}$, $ y(x_2) \in \{0.45,0.50,0.55,0.60\}$, $ y(x_3) \in \{0.20,0.21,0.22, \cdots 0.39\}$  and $ y(x_4) \in \{0.00,0.01,0.02,\cdots 0.19\}$ and $ y(x_5) = \pi $.  
The best path obtained is 
\begin{align*}
f_F (x) &=  \sum_{i=0}^{\myN} a_i x^i	,	
\end{align*}
with 
\begin{align*}
 a_0 &= 2,  \qquad & a_1 &= - 3.05046974259466, \\
 a_2 &=2.96364462153838, \qquad & a_3 &= - 1.7066436970901, \\
 a_4 &=0.481217918188328, \qquad & a_5 &=-0.0510588069227091,
\end{align*}
and the corresponding minimum time $ 1.00852712176 $ sec. We note that the error in the minimum time calculated is about $ 0.5\% $ and the total number of paths considered is $ 5000 $ for the random search and subsequently $ 3249 $  for the exhaustive search, i.e., a total of $8249$ paths.

\subsubsection{\textbf{Quantum exhaustive search algorithm}}\label{subsect_exhaustive_quantum}
As noted earlier, quantum search algorithms  give a quadratic speed-up compared to the classical algorithms in search problems. Therefore, the quantum search algorithm, Algorithm~$ \ref{Algo_2} $, could be employed to perform an exhaustive search.  We remark that the computational cost of Algorithm~$ \ref{Algo_2} $ in performing an exhaustive search on the states resulting from the same discretization as considered in Example A, is $ 4135 $ whereas the cost for the equivalent classical exhaustive search was $ 2825761 $ as noted earlier (see the Table~$ \ref{tab:title} $). 

\subsubsection{\textbf{Quantum random search algorithm}} \label{subsect_random_quantum}
The random search that we considered earlier using the classical algorithm in Example E, \mref{subsect_random}, can be carried out using the quantum search algorithm. If we simply pick $ 5000 $ random paths and apply the quantum search algorithm, Algorithm~$ \ref{Algo_2} $, on these paths to select the optimum path, then the cost of the quantum algorithm would be $ \sqrt{5000} \, \epsilon \approxeq 174 $ paths, where according to Theorem 6.2, \cite{bulger2003implementing}, the constant $ \epsilon $ has be taken to be approximately equal to $ 2.46 $

\subsubsection{\textbf{Quantum hybrid search algorithm}}\label{subsect_hybrid_quantum}
The hybrid quantum algorithm considered here is a two step quantum algorithm. The first step of this algorithm is an application of the quantum random search, in which the quantum algorithm, Algorithm~$ \ref{Algo_2} $, is applied on a predetermined number of states selected uniformly at random out of the total $ N$ states. The first step results in an approximate solution of the trajectory optimization problem. In the second step, we propose to employ Algorithm~$ \ref{Algo_2} $ again, following the random search,  for performing an exhaustive search in a finer discrete grid near the approximate solution obtained earlier. Clearly such a method combines the best features of both the randomized algorithm and the exhaustive quantum search algorithm and can be advantageous in many situations. If we continue with the example of the random quantum algorithm discussed in the previous subsection as the first step of the hybrid quantum algorithm then the corresponding cost of the first step would be $ \sqrt{5000} \, \epsilon \approxeq 174 $. The computational cost for the subsequent step (with application of Algorithm~$ \ref{Algo_2} $ again but on the paths considered in Example F, \mref{subsect_hybrid}) would be $ \sqrt{3249} \, \epsilon \approxeq 140  $. Hence, the total computational cost for the hybrid quantum algorithm is $ \approxeq 174 + 140 = 314 $.
A summary of the costs associated to various methods is given in the Table~$ \ref{tab:title} $.

\subsection{\textbf{Isoperimetric problem}}
As our second example, we consider the well-known isoperimetric problem \cite{bliss1925calculus}. The objective in this problem is to find the maximum area enclosed between a curve of the given fixed length and a given fixed straight line $ \ll $, such that the endpoints of the curve lie on the straight line $ \ll $. This problem could be solved by employing the principles of calculus of variations. Here we will obtain a numerical solution of this problem based on our proposed framework.

Let $ r = r(\theta) $, for $ \theta = \frac{\pi}{2}  $ to $ \theta = {\pi} $,  be a representation of the curve in polar coordinates with the boundary condition 
$ r(\frac{\pi}{2}) = 0 $. Let the length of the curve be $ \frac{\pi}{3} $, i.e., 
\begin{align}\label{length_condition}
	\int_{\pi/2}^{\pi} \sqrt{r^2 + \left(\frac{dr}{d \theta}\right)^2} \, d\theta = \frac{\pi}{3}.
\end{align}
Let the fixed straight line be the line $ \theta = \pi  $.  The objective function to be maximized is the area $ A $ given by
\begin{align}
	A = \int_{\pi/2}^{\pi} \, \frac{1}{2} r ^2 \, d\theta.
\end{align}

Next we explain our discretization scheme for this problem. The variable $ \theta  $ is treated as the independent variable and the variable  $r = r(\theta) $ is assumed to be the dependent variable. The parameter space consisting of the rectangle $ [\pi/2,\pi] \times [0,b] $ in $ \RR^2 $ is discretized, with a suitably chosen $ b \leq \frac{\pi}{3}$. We let $ \theta_k = \frac{k\pi}{\myN} $ for $k = 0,1 \cdots \myN $. Let $ r_{k} = r(\theta_k) \in \{ r_{[i]} = \frac{b i}{L} :\: i = 0,1 \cdots L \} $ for $ k=1 $ to $ \myN-1 $. We set $ r(\pi/2) = 0 $  to take into account the boundary condition. Further, using Lagrange interpolation we can determine the function $ r= r(\theta)$ as $r = c \sum_{k=0}^{\myN} r_{k} \phi_k(\theta) $ with $ c $ chosen such that \meqref{length_condition} is satisfied. 
Here $ \phi_k(\theta)  $  is chosen such that 
$ \phi_k(\theta_j) = 1 $ if $ j=k $ and $ \phi_k(\theta_j) = 0$ otherwise. More explicitly 
\begin{align}
	\phi_k(\theta) = \prod_{j=0,~ j\neq k}^{\myN} \, \frac{(\theta-\theta_j)}{\theta_k-\theta_j}.  
\end{align}
Essentially, it means that once a discretization of $ \theta $ is carried out then we require that at $\theta= \theta_k $, the corresponding $ r(\theta_k) $ can only have values from the set $ \{ r_{[i]} = \frac{b i}{L} :\: i = 0,1 \cdots L \} $. We note that with the above discretization there are $ N = (L+1)^{\myN} $  total possible paths. We fix $ b=1 $ and $ L= 100$ for the following examples.

\subsubsection{\textbf{Global discretization in physical space}} \label{sub_sect_discrete_physical_space_area_problem}
Examples G, H and I are given below to show the effect of increasing grid resolution on the accuracy of solutions. 

\begin{figure}
	\centering
	\includegraphics[scale=0.50]{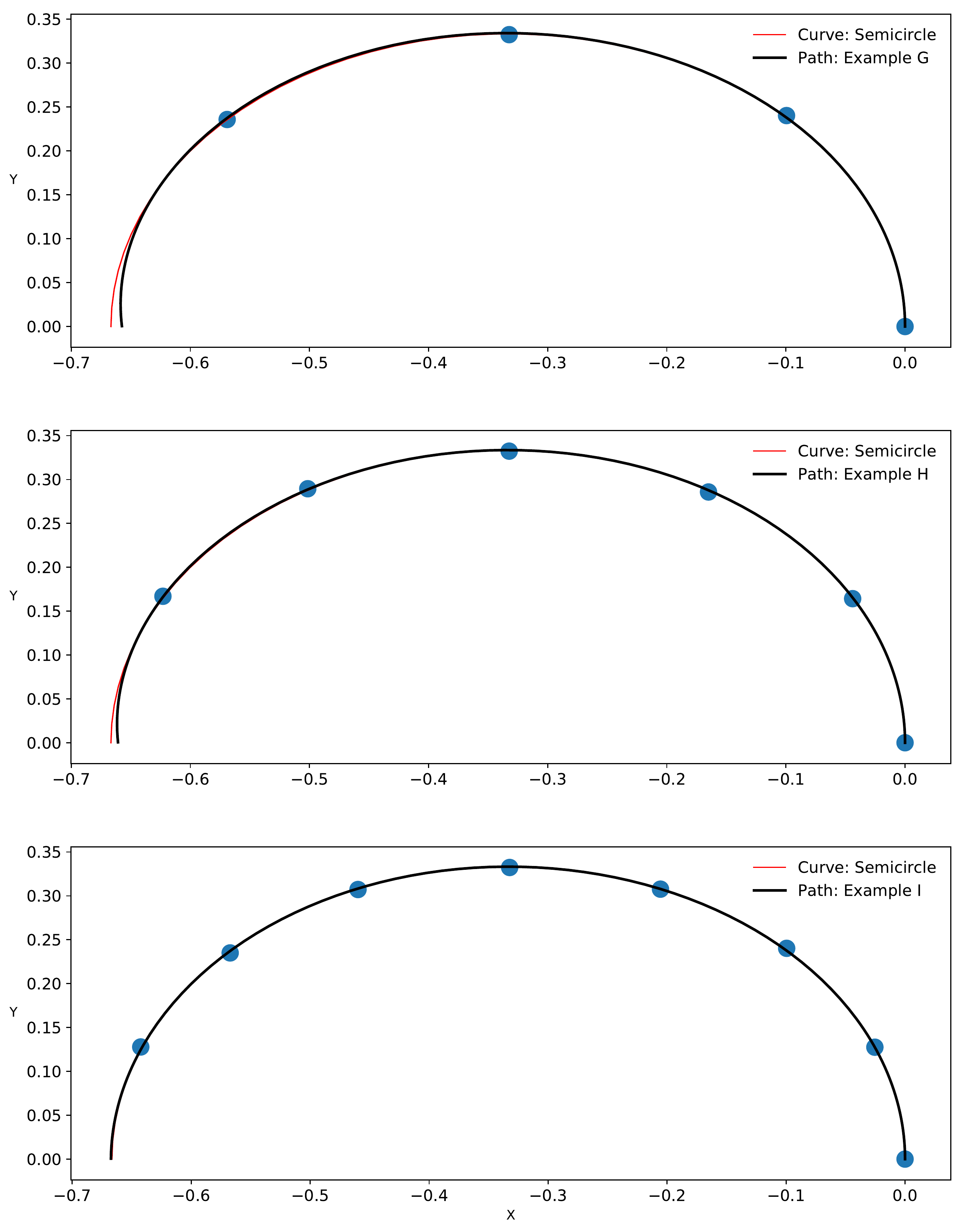}
	\caption{The paths for Examples G, H and I are plotted along with the semicircle resulting from the correct analytical solution of the isoperimetric problem. The nodes in the above examples represent the solution points on the discretization grid (see Table \ref{Table_Example_G}, Table \ref{Table_Example_H}, and Table \ref{Table_Example_I}).}\label{Fig_examples_G_H_I}
\end{figure}

\begin{table}
\begin{center}
	\small\setlength{\tabcolsep}{5pt}
	\begin{tabular}{c*{4}{c}}
		\toprule \toprule
		& \multicolumn{4}{c}{Example G} \\ 
		\hline
		\midrule
		$ \theta $:& $ \dfrac{\pi}{2} $  & $ \dfrac{5\pi}{8} $ &  $ \dfrac{3\pi}{4}$ & $ \dfrac{7\pi}{8} $ \\
		\midrule 
		$ r $: & $ 0$ & $ 0.26 $ &  $ 0.47$ &  $ 0.62$ \\
		\bottomrule
	\end{tabular} 
	\caption{Uniformly spaced $ 4 $ sample points in $ \theta $ for the isoperimetric problem.} \label{Table_Example_G}
\end{center}
\end {table}

\begin{table}
	\begin{center}
		\small\setlength{\tabcolsep}{5pt}
		\begin{tabular}{c*{6}{c}}
			\toprule \toprule
			& \multicolumn{6}{c}{Example H} \\ 
			\hline
			\midrule
			$ \theta $:& $ \dfrac{\pi}{2}   $  & $ \dfrac{7\pi}{12} $ & $ \dfrac{2\pi}{3} $ & $\dfrac{3\pi}{4} $ & $ \dfrac{5\pi}{6} $ & $ \dfrac{11\pi}{12} $ \\
			\midrule 
			$ r $: & $ 0 $ & $ 0.17 $ &  $ 0.33 $ &  $ 0.47 $ &  $ 0.58 $ & $ 0.64 $ \\
			\bottomrule
		\end{tabular} 
		\caption{Uniformly spaced $ 6 $ sample points in $ \theta $ for the isoperimetric problem.} \label{Table_Example_H}
	\end{center}
	\end {table}

	\begin{table}
		\begin{center}
			\small\setlength{\tabcolsep}{3pt}
			\begin{tabular}{c*{8}{c}}
				\toprule \toprule
				& \multicolumn{8}{c}{Example I} \\ 
				\hline
				\midrule
				$ \theta $:& $ \dfrac{\pi}{2} $ & $\dfrac{9\pi}{16} $ & $ \dfrac{5\pi}{8} $ & $ \dfrac{11\pi}{16} $ & $ \dfrac{3\pi}{4} $ & $ \dfrac{13\pi}{16} $ & $ \dfrac{7\pi}{8} $ &  $ \dfrac{15\pi}{16}  $ \\
				\midrule 
				$ r $: & $ 0 $ & $ 0.13 $ & $ 0.26 $ & $ 0.37 $ & $ 0.47 $ & $ 0.55 $ & $ 0.62 $ & $ 0.65 $ \\
				\bottomrule
			\end{tabular} 
			\caption{A finer grid resolution for obtaining a more accurate solution for the isoperimetric problem.} \label{Table_Example_I}
		\end{center}
		\end {table}

\begin{itemize}
	\item \textbf{Example G:} We consider the path obtained by using the Lagrange interpolation using the data of polar-coordinates $ (r,\theta) $ contained in Table \ref{Table_Example_G}, 
	\begin{align*}
	r_G (\theta) &= \pi \, \sum_{i=0}^{\myN} a_i  \theta^i	
	\end{align*}
	with 
	\begin{align*}
	a_0 &= -0.241651718296050,  \qquad & a_1 &= 0.0224680794862180, \\
	a_2 &=0.128732612832473, \qquad & a_3 &= - 0.0287105866008127, 
	\end{align*}

	This path can be discovered by considering $ \myN=3 $ in the preceding discussion. The size of the discrete search space  is $ 101^4  $. 
	We note that the maximum area obtained in this case, $ 0.174442279371647$ sq.~units, differs from the correct analytical solutions $ 0.17453292519 $ sq.~units by only about  $0.052$\%. \\
	
	\item \textbf{Example H:} Let us now consider a different discretization and pick the points as shown in Table \ref{Table_Example_H}. The path obtained by using the Lagrange interpolation in this case is  
	\begin{align*}
	r_G (\theta) &=  \pi \sum_{i=0}^{\myN} a_i \theta^i	,	
	\end{align*}
	with 
	\begin{align*}
	a_0 &= 0.117976693821150,  \qquad & a_1 &= - 0.776296300192081, \\
	a_2 &= 0.826825768983352, \qquad & a_3 &= - 0.328818088810480, \\
	a_4 &=0.0634200425140409, \qquad & a_5 &=- 0.00526623300376067.
	\end{align*}
This path can be discovered by using a discrete search space of size $ 101^6  $. The maximum area obtained is  $ 0.174481616034558 $ sq.~units. The error in the time obtained in this case is $ 0.029 $\%, which is better than the previous case.\\

	\item \textbf{Example I:}	 In this example, we get a more accurate result by considering even more points as shown in Table \ref{Table_Example_I}. The size of the discrete search space to discover this path, if we search exhaustively, is $ 101^{8} $. The maximum area obtained in this case is $0.174531915079274$ sq.~units and the error percentage is less than $ 0.0006 $\%.
\end{itemize}

\subsubsection{\textbf{Quantum exhaustive search algorithm}}\label{subsect_exhaustive_quantum_area}
The quantum search algorithm, Algorithm~$ \ref{Algo_2} $, could be employed to perform an exhaustive search using the same discretization as considered in Examples G, H and I to gain significant improvement in the size of the search space as discussed earlier.  In fact the size of the search space on using Algorithm~$ \ref{Algo_2} $, turns out to be $ 101^2 \, \epsilon $, $   101^3 \, \epsilon $   and $  101^4 \, \epsilon $, for the Examples G, H and I respectively, with $ \epsilon = 2.46 $ (see the Table~$ \ref{tab:title_area} $). 

\begin{table}
	\caption {A comparison of performances of different search methods for the  isoperimetric problem. (Assume $ \epsilon = 2.46 $.)} \label{tab:title_area} 
	\begin{center}
		\begin{tabular}{lllll}\toprule \toprule
			{\small Search }& {\small 	Search method description} & {\small Cost} & {\small Maximum area}  & {\small  Error} \% \\ 
			{\small method} & & & & \\
			\hline  \\
			{\small I} &	Classical exhaustive: Example G  & $ 101^4 $  & 0.174442279371647 &  0.052 \\ \\
			{\small I} &	Classical exhaustive: Example H   &   $ 101^6 $   &   0.174481616034558 &  0.029  \\ \\
			{\small I} &	Classical exhaustive: Example I   &   $ 101^8 $   &  0.174531915079274 & 0.0006 \\ \\
			{\small IV}&	Quantum exhaustive: Example G & $  101^2 \, \epsilon $ & 0.174442279371647 &0.052 \\ \\ 
			{\small IV}&	Quantum exhaustive: Example H & $  101^3 \, \epsilon $ &  0.174481616034558 & 0.029 \\ \\ 
			{\small IV}&	Quantum exhaustive: Example I & $  101^4 \, \epsilon $ & 0.174531915079274 & 0.0006 \\ \\ 
			\bottomrule
		\end{tabular} 
	\end{center}
\end{table}

\subsection{\textbf{Moon landing problem}} 

Our third example is the Moon landing problem considered in \cite{fahroo2002direct} and \cite{murecsan2012soft}. The problem is to find the optimum control for soft landing a spacecraft on the surface of the moon such that the fuel consumption is the minimum. We assume that the motion is vertical and that the lunar gravitation is a constant, $ g =1.63 $, throughout the motion. 
Let $ m(t)$ be the mass of the spacecraft including the fuel at the time $ t $. Let  $ h(t) $ and $ v(t) $ denote the height and the velocity of the spacecraft at the time  $ t $. The thrust at time $ t $, say $ T(t)$, is given by $ - k \frac{dm(t)}{dt} = -k u(t) $ and we treat $ u(t) = \frac{dm(t)}{dt} $ as a control variable.
The control problem can now be described as 
\begin{align}
\text{Maximize} \quad m(t),
\end{align}
subject to the equations of motion
\begin{align} \label{eq_motion}
	\frac{dh(t)}{dt}  = v(t), \qquad \frac{d^2h(t)}{dt^2} = -g - k \frac{u(t)}{m(t)},
\end{align}
and the constraints 
\begin{align}
m(t) \geq m_a = \text{the mass of the spacecraft without any fuel},\, -\mu \leq u(t)  \leq 0, 
\end{align}
where $ \mu $ is a constant which determines the maximum thrust.
Also, suppose the total time of flight is $\tau  $ then the condition for the soft landing is that  $ v(\tau) = 0 $ and $ h(\tau) = 0 $ and $ h(t) > 0  $ for $ t < \tau $. 
The initial condition is given by
\begin{align}
m(0) = m_0, \, v(0) = v_0 \, \text{and} \, h(0) = h_0. 
\end{align}
It easily follows on integrating the second equation of motion in  \meqref{eq_motion}  that $$ m(t) = m_0 \exp\left(\frac{v_0 -g\tau}{k}\right).$$
Therefore, instead of maximizing $ m(t) $ one can consider the equivalent objective function of minimizing the time of flight  $ \tau $.

	\begin{figure}[t]
		\centering
		\includegraphics[width=0.75\linewidth]{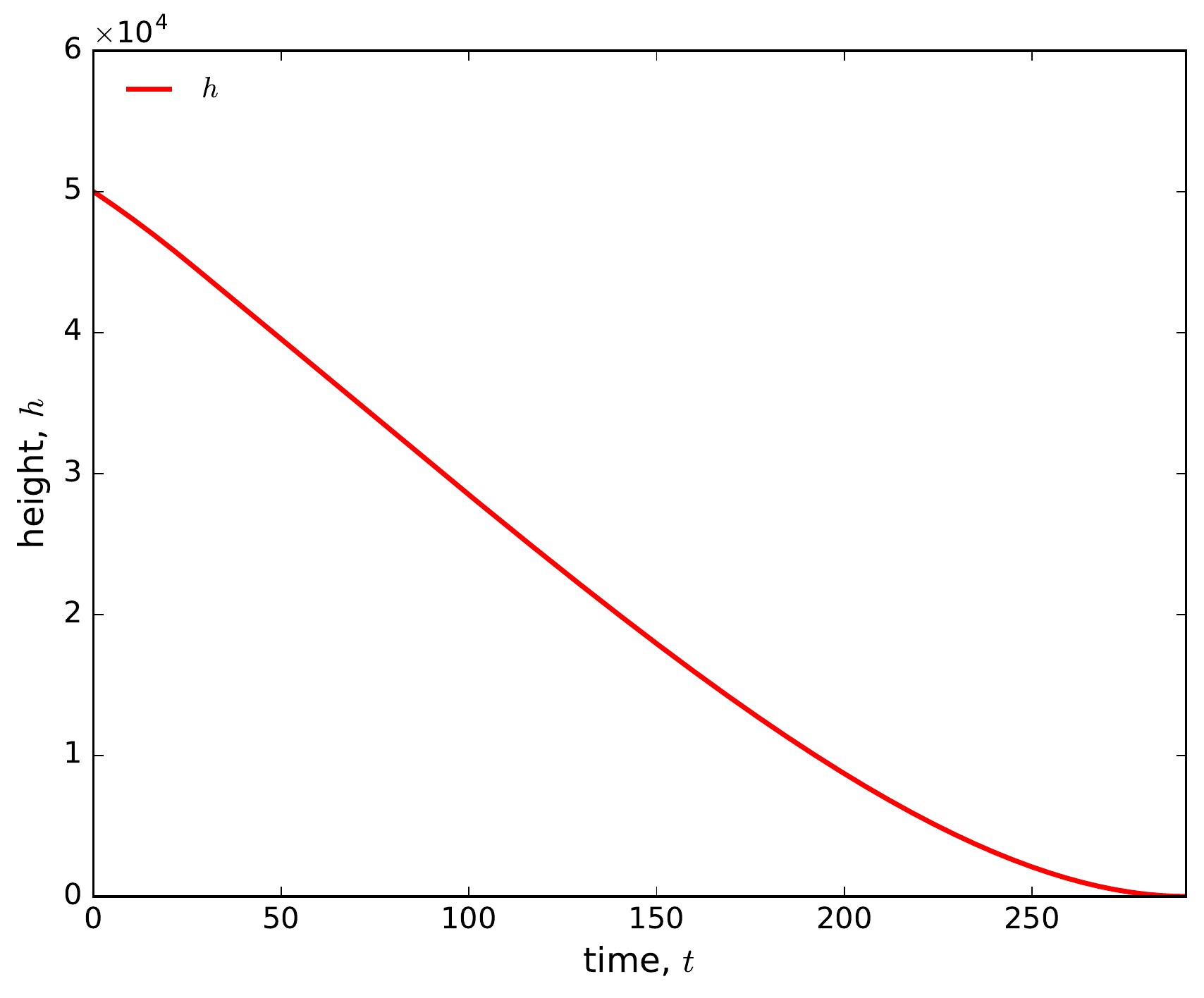} 
		\caption{{Time evolution of height for the moon landing problem}}\label{fig_h_t}
	\end{figure}
	\hfill

	\medskip
	
	\begin{figure}[t] 
				\centering
		\includegraphics[width=0.60\linewidth]{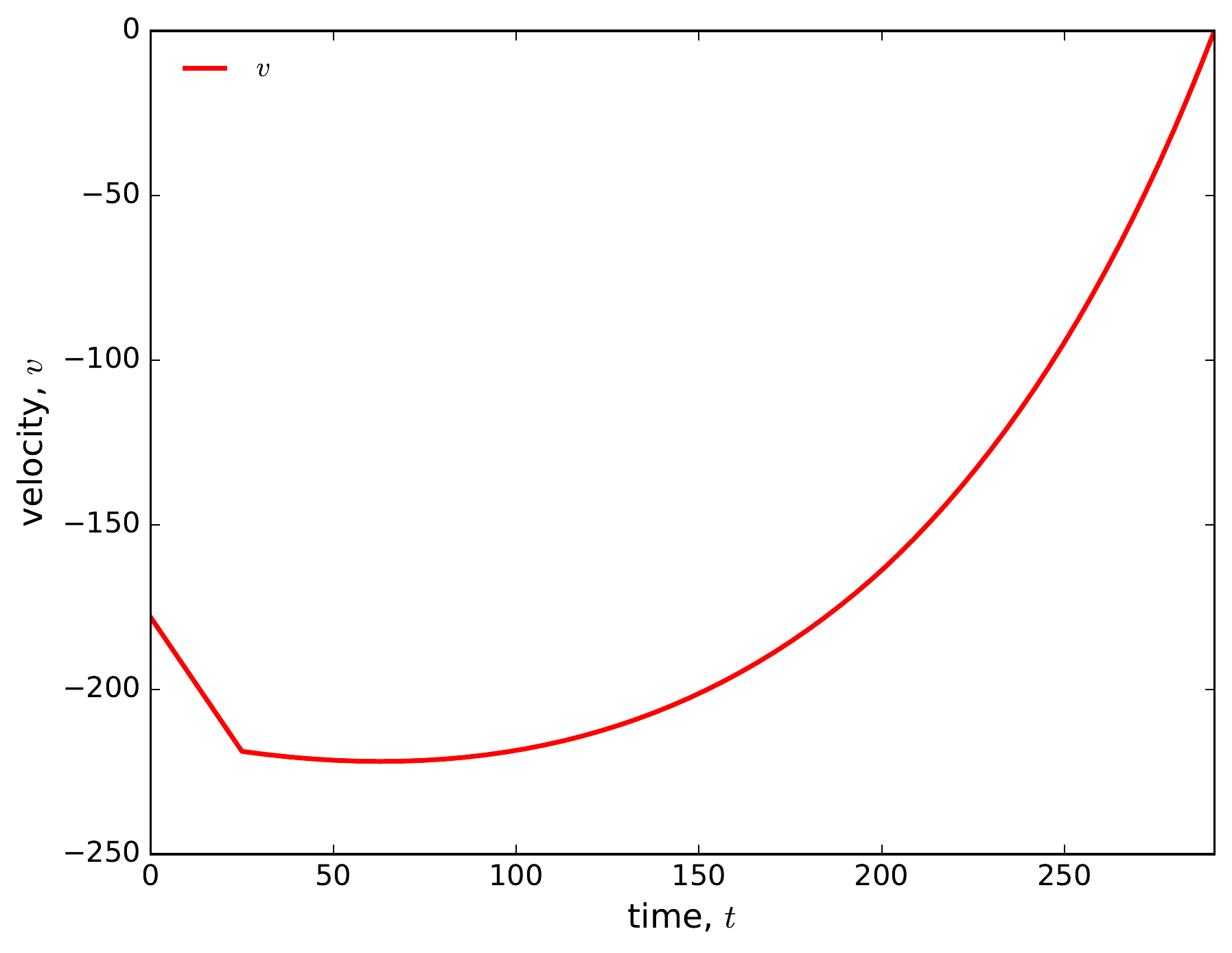} 
		\caption{{Time evolution of velocity for the moon landing problem (note that the straight line represents the free fall before the engine is switched on)} }\label{fig_v_t}
	\end{figure}

	\begin{figure}[t]
 		\centering
		\includegraphics[width=0.60\linewidth]{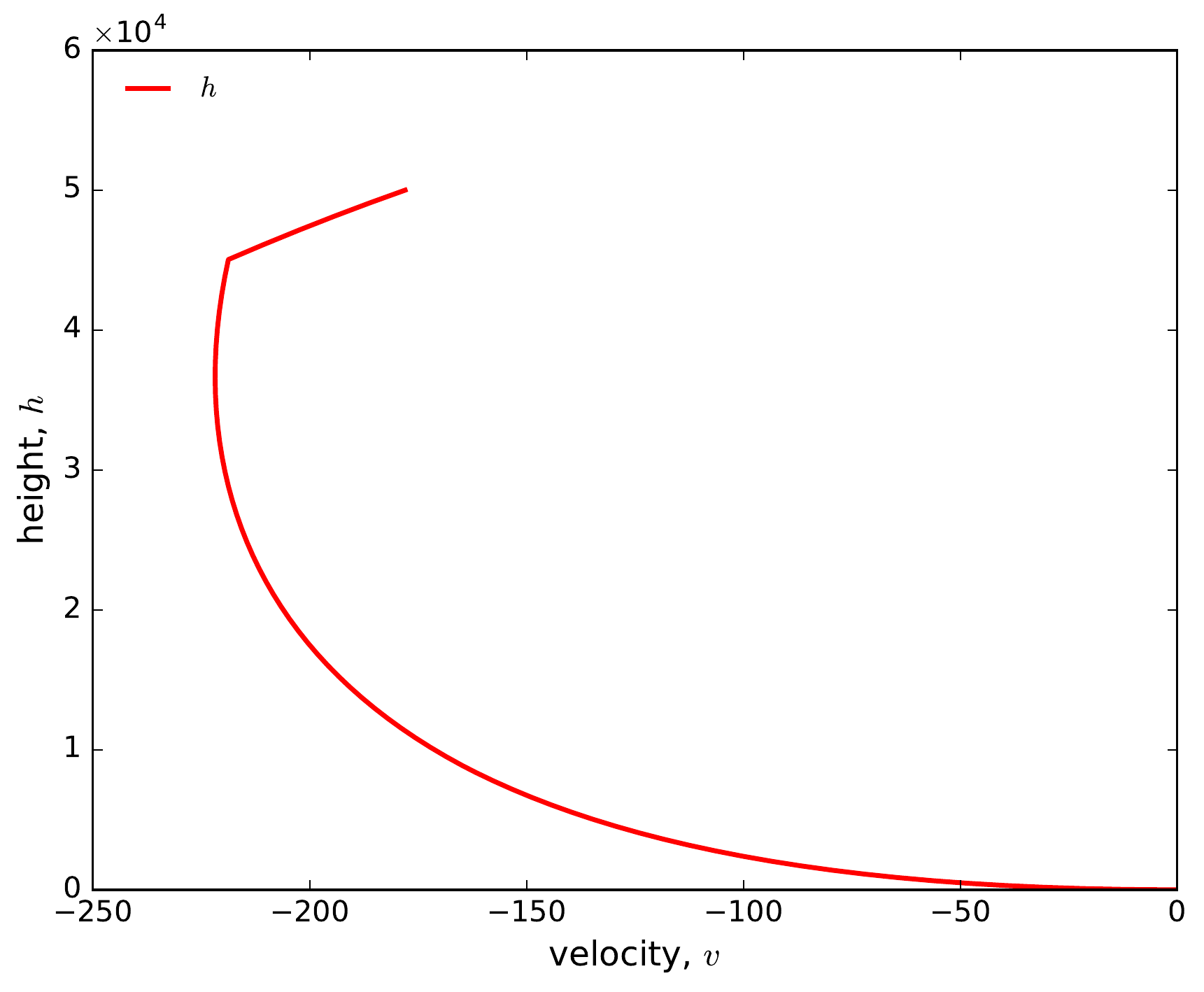} 
		\caption{{Variation of height and velocity for the moon landing problem}}\label{fig_h_v}
	\end{figure}
	
	\medskip
	\begin{figure}[t]	
		\centering
		\includegraphics[width=0.60\linewidth]{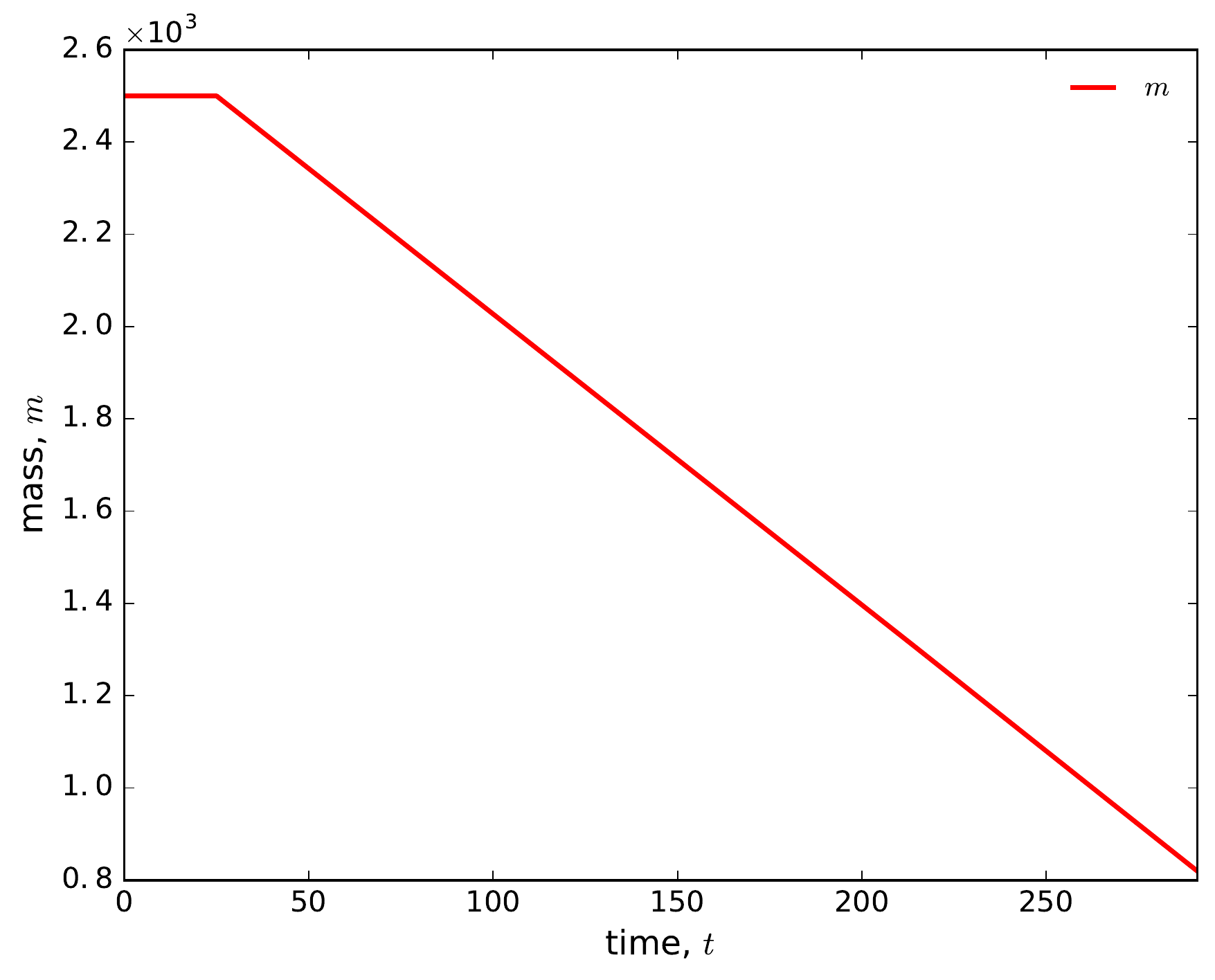} 
		\caption{{Time evolution of mass for the moon landing problem}}\label{fig_m_t}	
	\end{figure}

Next we describe our method for solving this equivalent control problem. We consider the following discretization. We pick $ 0 = t_0 < t_1 < t_2 < t_3 < t_4 = t_{\text{max}} $ as time points for discretization of the problem. We set 
\[ 
u(t) = \begin{cases} 
a_1 &  0 \leq t \leq t_1 \\
a_2  & t_1 < t  \leq t_2 \\
a_3  & t_2 < t \leq t_3 \\
a_4 &  t_3 < t \leq  t_4 
\end{cases}.
\]
Here, $ a_2, a_3$ and $ a_4 $ are chosen from a discrete set $ \{\alpha_k \} $ for $ k=1 $ to $ n $ with $ -\mu \leq \alpha_k < 0 $ and $ a_1 $ is to be determined later.
We note that, now with the help of equation of motions and the initial conditions $ m(t), v(t) $ and $ h(t) $ could be completely determined. In fact we get, 
\begin{align*}
m(t) = \begin{cases} 
a_k (t-t_{k-1}) + m(t_{k-1}), & \quad t_{k-1} \leq t  \leq  t_k \qquad \text{for } k=1,2,3,4.   
\end{cases}
\end{align*}
and 
\[ 
v(t) = \begin{cases} 
v_0 - g t -  k \log \left[\frac{a_k (t-t_{k-1}) + m(t_{k-1})}{m_0}\right], &  t_{k-1} \leq t  \leq  t_k  \quad \text{for } k=1,2,3,4.
\end{cases}
\]
Next, on integrating $ v(t) $ we get $ h(t) =  h_0 + \int_{0}^{t} v(x) \, dx$.  We treat $ a_1 $ and $ \tau $ as a variable to solve the constraint equations for a soft landing, i.e., we solve for $ \tau $ and $ a_1 $ such that $ h(\tau) = 0  $ and $ v(\tau) = 0 $ with $ h(t) > 0 $ for $ t < \tau $. We also let $ a_2,a_3 $ and $ a_4 $ take value in the discrete set $ \{-15,-14.99,-14.98,\cdots 0.1 \} $.  For the initial conditions $ h_0 = 50002.65, v_0 = -178, 
m_0 = 2500 $ and the sampling points $  t_1 = 25, t_2 = 75, t_3 = 200 $ and $ t_4 =400 $ with $ g = 1.63, k = 585, \mu = -15, m_a = 800 $ we obtain the following solution: 
$ a_1 = -0.0018, a_2 = -6.30, a_3 = -6.30, a_4 = -6.35 $ and $ \tau = 290.98 $ and $ m(\tau) = 819.76$. The time evolution of the mass, velocity and height as well as the variation of height versus velocity are shown in figures~\ref{fig_v_t} through \ref{fig_m_t}. These results  appear to be in agreement with the known solutions of this problem. We also note that these results are  based on Method I (see \mref{subsect_classical_exhaustive}). All the other methods (Method II to Method VI) may also be applied to solve this problem giving similar results. The performance comparisons of these methods are also expected to be similar those given in Table~$\ref{tab:title} $ for the brachistochrone problem, with quantum algorithms outperforming their classical counterparts.
\\ \\ \\

\noindent \textbf{Some remarks:} 
\begin{enumerate}
	\item We have arbitrarily chosen  $ a_1 $ as the variable chosen for satisfying the constraints for a soft landing. One can easily make an alternative choice.
	\item The choices of the time instants $ t_1, t_2, t_3 $   and $ t_4 $  were also arbitrary. In fact, one can let $ t_1,t_2,\cdots t_4 $ take values in an appropriately chosen  discrete set and then pick the best possible solutions. 
\end{enumerate}

\section{\textbf{Conclusion}} \label{Conclusion}
We have shown that a trajectory optimization problem (or a problem involving calculus of variations) can be formulated in the form of a search problem in a discrete space. An important feature of this work is our formulation of the discretization of the optimization problem incorporating the treatment of the dependent variable space in terms of appropriately chosen coarse grained states. The number of the coarse grained states is chosen to be significantly less than the traditional digital computer representations. This not only reduces the computational cost but also enables us to employ deterministic and stochastic approaches for obtaining global optimum. In particular, we presented the use of our proposed discretization approach and classical methods (Methods I--III) to solve the trajectory optimization problem. The proposed framework also enables the use of quantum computational algorithms for global (trajectory) optimization.  In this work, we showed that the discrete search problem can be solved by a variety of quantum algorithms including a quantum  exhaustive search in physical space or the coefficient space (Method IV), a quantum randomized search algorithm (Method V), or by employing a quantum hybrid algorithm (Method VI) depending on the nature of the problem. Quantum search algorithms offer a quadratic speed-up (in comparison to the traditional non-quantum approaches) and may well become methods of choice in many optimization applications once quantum computers become widely accessible and reach their true potential. 

\newpage

		\bibliographystyle{plain}
		\bibliography{Trajectory.bib}

\end{document}